%&amstex          
\input amstex\documentstyle{amsppt}  
\pagewidth{12.5cm}\pageheight{19cm}\magnification\magstep1
\topmatter
\title Graded Lie algebras and intersection cohomology\endtitle
\author G. Lusztig\endauthor
\address{Department of Mathematics, M.I.T., Cambridge, MA 02139}\endaddress
\thanks{Supported in part by the National Science Foundation.}\endthanks
\endtopmatter   
\document
\define\ufP{\un{\fP}}
\define\bcz{\ov{\cz}}
\define\occ{\ov{\cc}}

\define\uk{\un{\k}}

\define\tik{\ti{\k}}
\define\tcs{\ti{\cs}}

\define\du{\dot u}

\define\dk{\dot\k}

\define\bco{\bar{\co}}

\define\uX{\un X}

\define\bW{\bar W}

\define\hx{\hat x}

\define\uI{\un I}

\define\frl{\forall}

\define\si{\sim}
\define\wt{\widetilde}
\define\sqc{\sqcup}

\define\qua{\quad}

\define\tcl{\ti\cl}

\define\tfI{\ti{\fI}}

\define\bK{\bar K}

\define\lb{\linebreak}

\define\eSb{\endSb}

\define\op{\oplus}
   
\redefine\sp{\spadesuit}
\define\part{\partial}
\define\em{\emptyset}
\define\imp{\implies}
\define\ra{\rangle}
\define\n{\notin}
\define\iy{\infty}
\define\m{\mapsto}
\define\do{\dots}
\define\la{\langle}
\define\bsl{\backslash}

\define\sm{\smallmatrix}
\define\esm{\endsmallmatrix}
\define\sub{\subset}    
\define\bxt{\boxtimes}
\define\T{\times}
\define\ti{\tilde}
\define\nl{\newline}
\redefine\i{^{-1}}
\define\fra{\frac}
\define\un{\underline}
\define\ov{\overline}
\define\ot{\otimes}
\define\bbq{\bar{\QQ}_l}

\define\Ad{\text{\rm Ad}}

\define\Ind{\text{\rm Ind}}

\define\ind{\text{\rm ind}}

\define\a{\alpha}
\redefine\b{\beta}
\redefine\c{\chi}
\define\g{\gamma}
\redefine\d{\delta}
\define\e{\epsilon}
\define\et{\eta}
\define\io{\iota}
\redefine\o{\omega}
\define\p{\pi}

\define\r{\rho}
\define\s{\sigma}
\redefine\t{\tau}

\define\k{\kappa}
\redefine\l{\lambda}

\define\x{\xi}

\define\vp{\varpi}
\define\vt{\vartheta}

\redefine\G{\Gamma}
\redefine\D{\Delta}
\define\Om{\Omega}

\define\Th{\Theta}

\define\Ph{\Phi}

\define\kk{\bold k}

\define\BB{\bold B}
\define\CC{\bold C}

\define\II{\bold I}

\define\NN{\bold N}

\define\QQ{\bold Q}

\define\ZZ{\bold Z}

\define\ca{\Cal A}

\define\cc{\Cal C}
\define\cd{\Cal D}
\define\ce{\Cal E}
\define\cf{\Cal F}
\define\cg{\Cal G}
\define\ch{\Cal H}

\define\cj{\Cal J}
\define\ck{\Cal K}
\define\cl{\Cal L}
\define\cm{\Cal M}

\define\co{\Cal O}
\define\cp{\Cal P}

\define\car{\Cal R}
\define\cs{\Cal S}
\define\ct{\Cal T}
\define\cu{\Cal U}

\define\cz{\Cal Z}

\define\ff{\frak f}

\redefine\fI{\frak I}

\define\fP{\frak P}

\define\ta{\ti a}
\define\tb{\ti b}
\define\tc{\ti c}

\define\te{\ti e}
\define\tf{\ti f}
\define\tg{\ti g}

\define\tit{\ti t}

\define\tA{\ti A}

\define\tP{\ti P}

\define\tX{\ti X}

\define\bS{\bar S}

\define\tce{\ti\ce}
\define\bul{\bullet}

\define\dcl{\dot{\cl}}

\define\prq{\preceq}
\define\oL{\overset\circ\to L}

\define\bP{\un{P}}

\define\bQ{\un{Q}}
\define\dL{\dot{L}}
\define\sqq{\sqrt q}
\define\wInd{\wt{\Ind}}
\define\pt{\text{\rm point}}

\define\tOm{\ti{\Om}}

\define\Gi{G^\io}
\define\ucj{\un{\cj}}
\define\uY{\un Y}
\define\BBD{BBD}
\define\DE{D}
\define\CG{CG}
\define\CU{L1}
\define\CUII{L2}
\define\CUIII{L3}
\define\CB{L4}
\define\GR{L5}
\define\ADV{L6}
\define\ADVIII{L7}
\define\ADIX{L8}
\define\KA{K}
\define\VI{V}
\head Introduction\endhead
Let $(G,\io)$ be a pair consisting of a reductive connected algebraic group $G$ over $\CC$ 
and a homomorphism of algebraic groups $\io:\CC^*@>>>G$. The centralizer $\Gi$ of $\io(\CC^*)$
in $G$ acts naturally (with finitely many orbits) on the $n$-eigenspace $L_nG$ of 
$\Ad(\io(\CC^*))$ on the Lie algebra of $G$. (Here $n\in\ZZ-\{0\}$.) If $\co$ is a 
$\Gi$-orbit on $L_nG$ and $\cl$ is a $\Gi$-equivariant irreducible local system on $\co$ 
then the intersection cohomology complex $K=IC(\bco,\cl)$ is defined and we are interested
in the problem of computing, for any $\Gi$-orbit $\co'$ contained in $\bco$ and any 
$\Gi$-equivariant irreducible local system $\cl'$ on $\co'$, the multiplicity 
$m_{i;\cl,\cl'}$ of $\cl'$ in the local system obtained by restricting to $\co'$
the $i$-th cohomology sheaf of $K$. 

The main purpose of this paper is to give an algorithm to produce combinatorially a square 
matrix whose entries are polynomials with coefficients given by the 
multiplicities $m_{i;\cl,\cl'}$. Note that we do not have a purely combinatorial 
proof of the fact that the algorithm does not break down. We can only prove that by using 
geometry. But this will not prevent a computer from carrying out the algorithm.

The method of this paper relies very heavily on \cite{\GR} where many of the needed 
geometric results are proved. Note that in \cite{\GR} another purely algebraic description
of the multiplicities above was obtained, which however did not provide an algorithm for
computing them. 

While the existence of the algorithm above has an intrinsic interest, it also implies (by 
results in \cite{\CG}, \cite{\CUII, 10.7}) a solution of a problem in representation
theory, namely that of computing the multiplicities with which simple modules of an affine
Hecke algebra appear in a composition series of certain "standard modules".

At the same time, as a biproduct of the algorithm we find a way to compute the dimensions
of weight spaces in certain standard modules over an affine Hecke algebra. (See 4.6.)

In Section 1 we describe the algorithm. In Section 3 we show (based on the geometric
preliminaries in Section 2) that the algorithm in Section 1 is correct and it indeed leads
to the desired matrix of multiplicities. In Section 4 we define among other things a 
partial order on the set of isomorphism classes of irreducible $\Gi$-equivariant local 
systems on the various orbits in $L_nG$. In 4.7 we give a formulation of our results in 
terms of a canonical basis and two PBW-bases which generalizes the theory of canonical 
bases \cite{\CB} in the plus part of a quantized enveloping algebra of type $A_n$.

{\it Notation.} The cardinal of a finite set $S$ is denoted by $|S|$. 

Let $\ca=\ZZ[v,v\i]$ where $v$ is an indeterminate.

Let $\kk$ be an algebraically closed field of characteristic $p\ge0$. All algebraic 
varieties are assumed to be over $\kk$. 

We fix a prime number $l$ invertible in $\kk$. Let $\bbq$ be an algebraic closure of
the field of $l$-adic numbers. We will say "local system" instead of "$\bbq$-local system".
If $\cf$ is an irreducible local system (or its isomorphism class) over a subvariety $Y$ of
an algebraic variety $X$ we set $S_\cf=Y$.

For a connected affine algebraic group $H$ let $U_H$ be the unipotent radical of $H$, 
$\un{H}=H/U_H$, $LH$ the Lie algebra of $H$, $Z_H$ the centre of $H$.

\head 1. An algorithm\endhead
\subhead 1.1\endsubhead
Throughout this paper we assume that we are given a connected reductive algebraic group $G$
and a homomorphism of algebraic groups $\io:\kk^*@>>>G$. We assume that either $p=0$ or $p$
is sufficiently large (as in the last paragraph of \cite{\GR, 2.1(a)}). We set
$\Gi=\{g\in G;g\io(t)=\io(t)g\qua\frl t\in\kk^*\}$, a connected reductive subgroup of $G$.
We have $LG=\op_{n\in\ZZ}L_nG$ where 
$$L_nG=\{x\in LG;\Ad(\io(t))x=t^nx\qua\frl t\in\kk^*\}.$$ 
More generally, for any closed connected subgroup $H$ of $G$ that is normalized by 
$\io(\kk^*)$ we set $H^\io=H\cap\Gi$; we have $LH=\op_{n\in\ZZ}L_nH$ where 
$L_nH=LH\cap L_nG$. For $n\in\ZZ$, the adjoint action of $G$ on $LG$ restricts to an action
of $\Gi$ on $L_nG$.

\subhead 1.2\endsubhead
In the remainder of this paper we fix a subset $\D$ of $\ZZ$ consisting of two non-zero
elements whose sum is $0$. We assume that either $p=0$ or $\D\sub(-\iy,p)$.

We say that $(G,\io)$ is {\it rigid} if for some/any $n\in\D$ there exists a homomorphism 
of algebraic groups $\g:SL_2(\kk)@>>>G$ such that 
$\g\left(\sm t^n&0\\0&t^{-n}\esm\right)=\io(t^2)\mod Z_G$ for any $t\in\kk^*$. In this 
case, let $C_G^\io$ be the nilpotent $G$-orbit in $LG$ such that the corresponding
unipotent class in $G$ contains $\g(u)$ for any non-trivial unipotent element 
$u\in SL_2(\kk)$.

Let $\cp$ be the variety of parabolic subgroups of $G$. Let 
$\cp^\io=\{P\in\cp;\io(\kk^*)\sub P\}$. If $P\in\cp^\io$ then $\io$ gives rise to a 
homomorphism 

$\kk^*@>>>\bP$, $t\m(\text{image of $\io(t)$ under }P@>>>\bP)$,
\nl
denoted again by $\io$. Hence $L_n\bP$, $\bP^\io$ are well defined in terms of this $\io$ 
(we have $L_n\bP=L_nP/L_nU_P$).

\subhead 1.3\endsubhead
Let $\ct_G^{cu}$ (resp. $\ct_G^{pr}$) be the set of isomorphism classes of $G$-equivariant
irreducible local systems on some nilpotent orbit in $LG$ which are {\it cuspidal} (resp. 
{\it primitive}) in the sense of \cite{\CU, 2.2} (resp. \cite{\GR, 2.7}). 

We have $\ct_G^{cu}\sub\ct_G^{pr}$. The classification of local systems in $\ct_G^{pr}$ can
be deduced from the known classification of local systems in $\ct_G^{cu}$. For example, if
$G$ is simple of type $E_8$ then $\ct_G^{pr}$ consists of two objects: one is in 
$\ct_G^{cu}$ and one is $\bbq$ over the $G$-orbit $\{0\}$. Let 
$$\cj_G=\{(P,\ce);P\in\cp^\io,\ce\in\ct_{\bP}^{cu},(\bP,\io)\text{ is rigid },
C_{\bP}^\io=S_\ce\}.$$
Now $\Gi$ acts on $\cj_G$ by $g:(P,\ce)\m(gPg\i,\Ad(g)_!\ce)$. Let $\ucj_G$ be the set of 
orbits of this action. Let $K_G$ be the $\QQ(v)$-vector space with basis 
$(\II_{\cs})_{\cs\in\ucj_G}$.

For a connected affine algebraic group $\cg$ let $X_\cg$ be the variety of Borel subgroups
of $\cg$ and let $\text{\rm rk}(\cg)$ be the dimension of a maximal torus of $\cg$. We set 
$$e_\cg=\sum_j\dim H^{2j}(X_\cg,\bbq)v^{2j},\qua
\vt_\cg=(1-v^2)^{\text{\rm rk}(\cg)}e_\cg\in\ZZ[v^2].$$
If $\cf\in\ct_G^{pr}$ and $\cg$ is the connected centralizer in $G$ of some element in 
$S_\cf$, we set

$r_\cf=v^{-\dim X_\cg}e_\cg\in\ca$.

\subhead 1.4\endsubhead
Let $Q\in\cp^\io$. Associating to the $\bQ^\io$-orbit of $(P',\ce)\in\cj_{\bQ}$ the
$\Gi$-orbit of $(P,\ce)\in\cj_G$ (where $P$ is the inverse image of $P'$ under
$Q@>>>\bQ$ and $\bP',\bP$ are identified in the obvious way) defines a map
$a_Q^G:\ucj_{\bQ}@>>>\ucj_G$. We define a $\QQ(v)$-linear map 
$f_Q^G:K_{\bQ}@>>>K_G$ by $\II_{\cs'}\m\II_{a_Q^G(\cs')}$ for any $\cs'\in\ucj_{\bQ}$. 

\subhead 1.5\endsubhead
We define a map $\mu:\cj_G@>>>\ct_G^{pr}$ by $(P,\ce)\m\cf$ where $\cf$ is as follows. We 
choose a Levi $M$ of $P$ and we identify $M$ with $\bP$ in the obvious way. Then $\ce$ 
becomes a local system $\ce_M$ on a nilpotent $M$-orbit $D$ in $LM$. Let $C$ be the 
nilpotent $G$-orbit in $LG$ that contains $D$ and let $\cf$ be the unique $G$-equivariant 
local system on $C$ such that $\cf|_D=\ce_M$. (See \cite{\GR, 2.7}.) Then 
$\cf\in\ct_G^{pr}$ is clearly independent of the choice of $M$. 

For $\cf\in\ct_G^{pr}$ let $\cj_G^\cf=\mu\i(\cf)$. Then $\cj_G^\cf$ is $\Gi$-stable; let 
$\ucj_G^\cf$ be the set of $\Gi$-orbits on $\cj_G^\cf$. We have a partition 
$\ucj_G=\sqc_\cf\ucj_G^\cf$ and a direct sum decomposition
$$K_G=\op_\cf K_G^\cf$$
where $\cf$ runs over $\ct_G^{pr}$ and $K_G^\cf$ is the subspace of $K_G$ spanned by
$\{\II_{\cs};\cs\in\ucj_G^\cf\}$.

For $\cf\in\ct_G^{pr}$ let $Y^\cf$ be the set of all 
$((P,\ce),(P',\ce'))\in\cj_G^\cf\T\cj_G^{\cf^*}$ such that $P,P'$ have a common Levi. Now 
$\Gi$ acts diagonally on $Y^\cf$; let $\uY^\cf$ be the set of orbits. Define 
$\t:\uY^\cf@>>>\ZZ$ by
$$\Om\m\t(\Om)=\dim\fra{L_nU_{P'}+L_nU_P}{L_nU_{P'}\cap L_nU_P}
-\dim\fra{L_0U_{P'}+L_0U_P}{L_0U_{P'}\cap L_0U_P}$$
where $((P,\ce),(P',\ce'))$ is any element of the $\Gi$-orbit $\Om$ and $n\in\D$. (The fact
that this definition is independent of the choice of $n$ in $\D$ is seen as in
\cite{\GR, 16.3}.) 

\subhead 1.6\endsubhead
For $(P,\ce)\in\cj_G^\cf$ we choose a Levi $M$ of $P$ that contains $\io(\kk^*)$ and let 
$\tP$ be the unique parabolic subgroup with Levi $M$ such that $P\cap\tP=M$. We have 
$(\tP,\tce)\in\cj_G^\cf$ for a unique $\tce$. Although $\tP$ is not uniquely defined by 
$P$, its $\Gi$-orbit is uniquely defined by the $\Gi$-orbit of $(P,\ce)$ (since $M$ is 
uniquely defined by $P$ up to the conjugation action of $U_P^\io$). Thus 
$(P,\ce)\m(\tP,\tce)$ induces a well defined involution $\cs\m\tcs$ of $\cj_G^\cf$. 

Similarly, for $((P,\ce),(P',\ce'))\in Y^\cf$ we choose a common Levi $M$ of $P$ and $P'$ 
that contains $\io(\kk^*)$ and let $\tP$ be the unique parabolic subgroup with Levi $M$ 
such that $P\cap\tP=M$. We have $((\tP,\tce),(P',\ce'))\in Y^\cf$ for a unique $\tce$. 
Again the $\Gi$-orbit of $((\tP,\tce),(P',\ce'))$ is uniquely defined by the $\Gi$-orbit of
$((P,\ce),(P',\ce'))$ (since $M$ is uniquely defined by $P,P'$ up to the conjugation action
of $U_{P\cap P'}^\io$). Thus $((P,\ce),(P',\ce'))\m((\tP,\tce),(P',\ce'))$ induces a well 
defined involution $\Om\m\tOm$ of $\uY^\cf$. As in the proof of \cite{\GR, 16.4(c)} we have
$$\t(\Om)+\t(\tOm)=c_\cf\tag a$$
where
$$c_\cf=\dim L_nG-\dim L_0G-\dim L_n\bP+\dim L_0\bP$$
for any $(P,\ce)\in\cj_G^\cf$. 
(If $(P,\ce')$ is another pair in $\cj_G^\cf$ then there exists an isomorphism
$\bP@>\si>>\bP'$ which is compatible with $\io$ so that $c_\cf$ depends only on $\cf$).

\subhead 1.7\endsubhead
Let $\cf\in\ct_G^{pr}$. In (a), (b) below we give a "combinatorial" interpretation of the
sets $\ucj_G^\cf$ and $\uY^\cf$. We may assume that $S_\cf\cap L_nG\ne\em$; otherwise
both our sets are empty.

Let $M$ be the centralizer in $G$ of some maximal torus of the connected centralizer in 
$\Gi$ of some element in $S_\cf\cap L_nG$. Then $\io(\kk^*)\sub M$ and $M$ is independent
of the choices (up to $\Gi$-conjugacy) since, by \cite{\GR, 14.5}, $S_\cf\cap L_nG$ is a
single $\Gi$-orbit. Let $\uX=\{P\in\cp;M\text{ is a Levi of }P\}$. 
If $P\in\uX$ then $(P,\ce)\in\cj_G^\cf$ for a unique $\ce$, see \cite{\GR, 11.6(c)}. We
have an imbedding $\uX@>>>\cj_G^\cf$, $P\m(P,\ce)$ and an imbedding 
$\uX@>>>\cj_G^{\cf^*}$, $P\m(P,\ce^*)$. Let $N_GM$ be the normalizer of $M$ in $G$. It
is known that the conjugation action of $N_GM/M$ on $\uX$ is simply transitive. Note that 
$(N_GM)^\io/M^\io$ is naturally a subgroup of $N_GM/M$. We show:

(a) {\it the map $j_1:\text{(set of $(N_GM)^\io/M^\io$-orbits on }\uX)@>>>\ucj_G^\cf$
induced by the imbedding $\uX@>>>\cj_G^\cf$ is bijective;}

(b) {\it the map $j_2:\text{(set of $(N_GM)^\io/M^\io$-orbits on }\uX\T\uX)@>>>\uY^\cf$
(orbits for diagonal action) induced by the imbedding $\uX\T\uX@>>>Y^\cf$ is bijective.}
\nl
Let $(P,\ce)\in\cj_G^\cf$. We can find a Levi $M'$ of $P$ that contains $\io(\kk^*)$. 
Using \cite{\GR, 11.4} we see that there exists $g\in\Gi$ such that $gM'g\i=M$. Then 
$gPg\i\in\uX$. We see that $j_1$ is surjective. 

Now let $P,P'\in\uX$ be such that $P'=gPg\i$ for some $g\in\Gi$. Then $M$ and $M'=g\i Mg$ 
are Levi subgroups of $P$ that contain $\io(\kk^*)$. There is a unique $u\in U_P$ such that
$uMu\i=M'$. For any $t\in\kk^*$ we set $u'=\io(t)u\io(t)\i\in U_P$; we have

$M'=\io(t)M'\io(t)\i=u'(\io(t)M\io(t)\i)u'{}\i=u'Mu'{}\i$.
\nl
By the uniqueness of $u$ we have $u'=u$. Thus $u\in U_P^\io$. Let $g'=gu$. Then $g'\in\Gi$,
$P'=g'Pg'{}\i$, $M=g'{}\i Mg'$. Thus, $g'\in(N_GM)^\io$. We see that $j_1$ is injective. 
This proves (a).

Let $((P,\ce),(P',\ce'))\in Y^\cf$. We can find a common Levi $M'$ of $P,P'$ that contains
$\io(\kk^*)$. As in the proof of (a) we can find $g\in\Gi$ such that $gM'g\i=M$. Then 
$(gPg\i,gP'g\i)\in\uX\T\uX$. We see that $j_2$ is surjective. 

Now let $P,P_1,P',P'_1$ in $\uX$ be such that $P'=gPg\i,P'_1=gP_1g\i$ for some $g\in\Gi$. 
Then $M$ and $M'=g\i Mg$ are Levi subgroups of $P\cap P_1$ that contain $\io(\kk^*)$. There
is a unique $u\in U_{P\cap P_1}$ such that $uMu\i=M'$. As in the proof of (a) we see using
the uniqueness of $u$ that $u\in U_{P\cap P_1}^\io$. Let $g'=gu$. Then $g'\in\Gi$, 
$P'=g'Pg'{}\i$, $P'_1=g'P_1g'{}\i$, $M=g'{}\i Mg'$. Thus, $g'\in(N_GM)^\io$. We see that 
$j_2$ is injective. This proves (b).

\subhead 1.8\endsubhead
Define a symmetric $\QQ(v)$-bilinear form $(:):K_G\T K_G@>>>\QQ(v)$ by setting 
(for $\cs\in\ucj_G^\cf$, $\cs'\in\ucj_G^{\cf'}$): 
$$(\II_{\cs}:\II_{\cs'})=0\text{ if }\cf'\ne\cf^*,$$
$$(\II_{\cs}:\II_{\cs'})=\fra{\vt_{\Gi}}{\vt_{Z_{\bP}^0}}
\sum\Sb\Om\in\uY^\cf\\ \p_1(\Om)=\cs,\p_2(\Om)=\cs'\endSb
(-v)^{\t(\Om)}\text{ if }\cf'=\cf^*;$$
here $\p_1:\uY^\cf@>>>\ucj_G^\cf$, $\p_2:\uY^\cf@>>>\ucj_G^{\cf'}$ are the obvious 
projections and $P\in\cp$ is such that $(P,\ce)\in\cs$ for some $\ce$.

Let $\bar{}:\QQ(v)@>>>\QQ(v)$ be the $\QQ$-algebra involution such that $\ov{v^m}=v^{-m}$ 
for all $m\in\ZZ$. Define a $\QQ$-linear involution $\b:K_G@>>>K_G$ by
$$\b(\r\II_{\cs})=\ov{\r}\II_{\cs}$$
for all $\r\in\QQ(v),\cs\in\ucj_G$. Define a $\QQ(v)$-linear involution 
$\s:K_G^\cf@>>>K_G^\cf$ by $\s(\II_{\cs})=\II_{\tcs}$ (see 1.6) for all $\cs\in\ucj_G^\cf$.
From 1.6(a) we see that, for $\x\in K_G^\cf$, $\x'\in K_G^{\cf'}$, we have
$$\ov{(\b(\x):\b(\x'))}=(-v)^{c_\cf}(\s(\x):\x').\tag a$$
Let 
$$\car_G=\{\x\in K_G;(\x:K_G)=0\}.$$ 
From (a) we see that $\b(\car_G)\sub\car_G$. Clearly,
if $Q\in\cp^\io$ then for $\x\in K_{\bQ}$ we have
$$f_Q^G(\b(\x))=\b(f_Q^G(\x)).\tag b$$
Moreover,
$$f_Q^G(\car_{\bQ})\sub\car_G.\tag c$$
See 3.5 for a proof.

\subhead 1.9\endsubhead
Let $n\in\D$. Let $P\in\cp^\io$ be such that $(\bP,\io)$ is rigid. We can find a Levi
subgroup $M$ of $P$ such that $\io(\kk^*)\sub M$. Let $s$ be the unique element in 
$[LM,LM]$ such that $[s,x]=mx$ for any $m\in\ZZ,x\in L_mM$. We have
$LG=\op_{r\in(n/2)\ZZ}L^rG$ where $L^rG=\{x\in LG;[s,x]=rx\}$ and
$LG=\op_{r\in(n/2)\ZZ,t\in\ZZ}L^r_tG$ where $L^r_tG=L^rG\cap L_tG$. We say that $P$ (as 
above) is $n$-good if 

$LU_P=\op_{r\in(n/2)\ZZ,t\in\ZZ;2t/n<2r/n}L^r_tG$.
\nl
(This implies that 

$LM=\op_{r\in(n/2)\ZZ,t\in\ZZ;2t/n=2r/n}L^r_tG$,
$LP=\op_{r\in(n/2)\ZZ,t\in\ZZ;2t/n\le2r/n}L^r_tG$.)
\nl
Note that the condition that $P$ is $n$-good is independent of the choice of $M$.

Let $\fP_n$ be the set of all $P\in\cp^\io$ such that $(\bP,\io)$ is rigid and $P$ is 
$n$-good. Let $\fP'_n=\{P\in\fP_n;P\ne G\}$. Now $\Gi$ acts on $\fP_n,\fP'_n$ by 
conjugation. Let $\ufP_n,\ufP'_n$ be the sets of orbits of these actions. These are finite
sets since $\Gi$ acts with finitely many orbits on $\cp^\io$. We have $G\in\fP_n$ if and 
only if $(G,\io)$ is rigid. Hence $\ufP'_n=\ufP_n$ if $(G,\io)$ is not rigid, 
$\ufP_n=\ufP'_n\sqc\{G\}$ if $(G,\io)$ is rigid.

\subhead 1.10\endsubhead
Let $n\in\D$. For $\et\in\ufP_n$ we set
$$d_\et=\dim L_0G-\dim L_0P+\dim L_nP$$
where $P\in\et$. For $\et,\et'\in\ufP'_n$ we say that $\et'\prec\et$ if $d_{\et'}<d_\et$. 
We say that $\et'\prq\et$ if either $\et=\et'$ or $\et'\prec\et$. Now $\prq$ is a partial 
order on $\ufP'_n$.

\subhead 1.11\endsubhead
Our goal is to define subsets $\cz^\et_n$ of $K_G$ (for $n\in\D$ and $\et\in\ufP_n$)). The
definition of these subsets is inductive and is based on a number of lemmas which will be 
verified in Section 3 (where we assume, as we may, that $\kk$ is an algebraic closure of a
finite field). If $G$ is a torus then $\et$ must be $\{G\}$ 
and $\cz^\et_n$ consists of the unique basis element of $K_G$. We now assume that $G$ is 
not a torus and that the subsets $\cz^\et_n$ are already defined when $G$ is replaced by 
any $\bP$ with $P\in\fP'_n$.

Our definition is based on the following scheme.

(i) We first define $\cz^\et_n$ in the case where $\et\in\ufP'_n$ by
$$\cz^\et_n=f_P^G(\cz^{\{\bP\}}_n)$$
 where $P\in\et$ and $f_P^G$ is as in 1.4. (Note that 
$\cz^{\{\bP\}}_n\sub K_{\bP}$ is defined by the inductive assumption.) We set 
$$\cz'_n=\cup_{\et\in\ufP'_n}\cz^\et_n.$$

(ii) Using (i) we define elements $W^\x_n$ for any $\x\in\cz'_n$ by a procedure similar to
the definition of the "new" basis of a Hecke algebra. (See 1.13.)

(iii) Using (i) and (ii) for $n$ and $-n$ we define $\cz^\et_n$ for $\et\in\ufP_n-\ufP'_n$.
(See 1.18.)

\proclaim{Lemma 1.12} Let $n\in\D$. 

(a) If $\et\in\ufP'_n$ and $P\in\et$, the map $\cz^{\{\bP\}}_n@>>>\cz^\et_n$ given by 
$\x\m f_P^G(\x)$ is bijective.

(b) The union $\cup_{\et\in\ufP'_n}\cz^\et_n$ is disjoint.

(c) In the setup of (a) let $\x'\in\cz'_n$ (relative to $\bP$ instead of $G$). Then 
$f_P^G(\x')$ is an $\ca$-linear combination of elements in various $\cz^{\et'}_n$ (with
$\et'\in\ufP'_n$, $\et'\prec\et$) plus an element of $\car_G$.

(d) In the setup of (a) let $\x_0\in\cz^{\{\bP\}}_n$. Then $\b(\x_0)-\x_0$ is an 
$\ca$-linear 
combination of elements in $\cz'_n$ (relative to $\bP$) plus an element of $\car_{\bP}$.

(e) The matrix with entries $(\x:\x')\in\QQ(v)$ indexed by $\cz'_n\T\cz'_n$ is 
non-singular.
\endproclaim
See 3.6, 3.7 for a proof.

\subhead 1.13\endsubhead
Let $n\in\D$. We show that for any $\x\in\cz'_n$ we have
$$\b(\x)=\sum_{\x_1\in\cz'_n}a_{\x,\x_1}\x_1\mod\car_G\tag a$$
where $a_{\x,\x_1}\in\ca$ are uniquely determined and satisfy the following conditions
(where $\et,\et_1$ are given by $\x\in\cz^\et_n,\x_1\in\cz^{\et_1}_n$):

$a_{\x,\x_1}\ne0$ implies $\et_1\prec\et$ or $\x=\x_1$;

$a_{\x,\x_1}=1$ if $\x=\x_1$.
\nl
We have $\x=f_P^G(\x_0)$ where $\x_0\in\cz_n^{\{\bP\}}$ (notation of 1.11(i)). We express 
$\b(\x_0)-\x_0$ as in 1.12(d). Applying $f_P^G$ and using 1.12(c), 1.8(d) and 1.8(c) we 
deduce that (a) holds except perhaps for the uniqueness statement. To show the uniqueness 
we note that the $a_{\x,\x_1}$ are determined from the system of linear equations

$(\b(\x):\x_2)=\sum_{\x_1\in\cz'_n}(\x_1:\x_2)a_{\x,\x_1}$
\nl
(with $\x_2\in\cz'_n$) whose matrix of coefficients is invertible by 1.12(e).

Using the equality $\b^2=1:K_G@>>>K_G$ and the inclusion $\b(\car_G)\sub\car_G$ (see 1.8)
we see that for any $\x,\x_1$ in $\cz'_n$ we have

$\sum_{\x_2\in\cz'_n}\ov{a_{\x,\x_2}}a_{\x_2,\x_1}=\d_{\x,\x_1}$.
\nl
Using a standard argument we see that there is a unique family of elements 
$c_{\x,\x_1}\in\ZZ[v]$ (defined for $\x,\x_1\in\cz'_n$) such that for any 
$\x,\x_1\in\cz'_n$ with $\x\in\cz^\et_n,\x_1\in\cz^{\et_1}_n$ we have

$c_{\x,\x_1}=\sum_{\x_2\in\cz'_n}\ov{c_{\x,\x_2}}a_{\x_2,\x_1}$;

$c_{\x,\x_1}\ne0$ implies $\et_1\prec\et$ or $\x=\x_1$;

$c_{\x,\x_1}\ne0$, $\x\ne\x_1$ implies $c_{\x,\x_1}\in v\ZZ[v]$;

$c_{\x,\x_1}=1$ if $\x=\x_1$.
\nl
For $\x\in\cz'_n$ we set $W_n^\x=\sum_{\x_1\in\cz'_n}c_{\x,\x_1}\x_1$. Then 
$\b(W_n^\x)=W_n^\x\mod\car_G$.

\subhead 1.14\endsubhead
Until the end of 1.16 we assume that $(G,\io)$ is rigid. For any $\cf\in\ct_G^{pr}$ such 
that $S_\cf=C_G^\io$ we set $\cs_\cf=\mu\i(\cf)$ which, by \cite{\GR, 11.9}, is a single 
$\Gi$-orbit $\cj_G$ ($\mu$ as in 1.5). Let 
$$\cc'=\{r_\cf\i\II_{\cs_\cf};\cf\in\ct_G^{pr},S_\cf=C_G^\io\}\sub K_G;$$
here $r_\cf\in\ca$ is as in 1.3.

\subhead 1.15\endsubhead
Let $n\in\D$. For any $x\in K_G$ we define $Y_n(x)\in K_G$ by the conditions
$$(Y_n(x):\cz'_n)=0,\qua x=Y_n(x)+\sum_{\x\in\cz'_n}\g_\x\x$$
with $\g_x\in\QQ(v)$. The coefficients $\g_\x$ are determined from the system of linear 
equations

$\sum_{\x\in\cz'_n}(\x:\x')\g_\x=(x:\x')$
\nl
with $\x'\in\cz'_n$, whose matrix of coefficients is invertible by 1.12(e). Let 
$$J_{-n}=\{\x_0\in\cz'_{-n};Y_n(W_{-n}^{\x_0})\n\car_G\}$$
and let $\cc_n$ be the image of the map $J_{-n}@>>>K_G$, $\x_0\m Y_n(W_{-n}^{\x_0})$. This
can be regarded as a surjective map $h_n:J_{-n}@>>>\cc_n$.

\proclaim{Lemma 1.16} $h_n$ is bijective. 
\endproclaim
See 3.10 for a proof.

\proclaim{Lemma 1.17} For $n\in\D$, the union $\cz'_n\cup\cc_n\cup\cc'$ is disjoint.
\endproclaim
See 3.11 for a proof.

\subhead 1.18\endsubhead
If $(G,\io)$ is not rigid then $\ufP_n=\ufP'_n$ and the definition of the subsets 
$\cz_n^\et$ ($\et\in\ufP_n$) is complete. If $(G,\io)$ is rigid and $n\in\D$ we set 
$\cz_n^{\{G\}}=\cc_n\cup\cc'$. By 1.17, this union is disjoint. 
The definition of the subsets $\cz_n^\et$ ($\et\in\ufP_n$) is complete.

We set $\cz_n=\cz'_n$ if $(G,\io)$ is not rigid and $\cz_n=\cz'_n\cup\cz_n^{\{G\}}$ if 
$(G,\io)$ is rigid. By 1.17, the last union is disjoint. 

\subhead 1.19\endsubhead
For $n\in\D$ and $\x\in\cz_n$ we define an element $W_n^\x$ as follows. When $\x\in\cz'_n$,
this is already defined in 1.13. When $(G,\io)$ is rigid we set 
$W_n^\x=W_{-n}^{h_n\i(\x)}$ if $\x\in\cc_n$ and $W_n^\x=\x$ if $\x\in\cc'$. 

We now define a matrix $(c_{\x,\x'})$ with entries in $\QQ(v)$ indexed by $\cz_n\T\cz_n$ by
the following requirements:

When $\x,\x'\in\cz'_n$ then $c_{\x,\x'}$ are as in 1.13. 

When $\x\in\cz'_n,\x'\n\cz'_n$ then $c_{\x,\x'}=0$. 

When $\x\n\cz'_n,\x'\n\cz'_n$ then $c_{\x,\x'}=\d_{\x,\x'}$. 

When $\x\n\cz'_n$ then $c_{\x,\x'}$ for $\x'\in\cz'_n$ are determined by the system of 
linear equations $(W_n^\x:\x'')=\sum_{\x'\in\cz'_n}(\x':\x'')c_{\x,\x'}$ (with 
$\x''\in\cz'_n$) whose matrix of coefficients has invertible determinant.

Note that for any $\x\in\cz_n$ we have 
$$W_n^\x=\sum_{\x'\in\cz_n}c_{\x,\x'}\x'.\tag a$$

\proclaim{Lemma 1.20} Let $n\in\D$. Let $\cs\in\ucj_G$. There exist $e_{\cs,\x}\in\ca$ (for
$\x\in\cz_n$) and $r\in\car_G$ such that $\II_\cs=\sum_{\x\in\cz_n}e_{\cs,\x}\x+r$. 
\endproclaim
See 3.14 for a proof.

\head 2. Geometric preliminaries\endhead
\subhead 2.1\endsubhead
In this section we assume that $\kk$ is an algebraic closure of the finite field $F_p$ with
$|F_p|=p$. For $q\in\{p,p^2,\do\}$ let $F_q$ be the subfield of $\kk$ with $|F_q|=q$. If 
$X$ is an algebraic variety we denote by $\cd(X)$ the bounded derived category of 
(constructible) $\bbq$-sheaves on $X$. For $K\in\cd(X)$ let $\ch^iK$ be the $i$-th 
cohomology sheaf of $K$. For $n\in\ZZ$ let $\bbq(n/2)$ be as in the Introduction to
\cite{\ADIX}. We write $K[[n/2]]$ instead of $K[n]\ot\bbq(n/2)$. We fix a square root 
$\sqrt{p}$ of $p$ in $\bbq$. If $q$ is a power $p^e$ of $q$ we set $\sqrt{q}=(\sqrt{p})^e$.
We shall assume that Frobenius relative to $F_q$ acts on $\bbq(n/2)$ as multiplication by 
$(\sqrt{q})^{-n}$.

For a connected affine algebraic group $\cg$ we have 
$$\vt_\cg|_{v=-1/\sqq}=|L\cg(F_q)|\i|\cg(F_q)|$$
for any $F_q$-rational structure on $\cg$ such that $\un{\cg}$ is $F_q$-split; here 
$\vt_\cg$ is as in 1.3. 

Define $\o:\kk@>>>\kk$ by $x\m x^p-x$. Let $\cu$ be a local system of rank $1$ on $\kk$ 
such that $\cu\op\bbq$ is a direct summand of $\o_!\bbq$. Let $E,E'$ be two $\kk$-vector 
spaces of the same dimension $<\iy$ and let $\s:E\T E'@>>>\kk$ be a perfect bilinear 
pairing. Let $s:E\T E'@>>>E$, $s':E\T E'@>>>E'$ be the projections. Recall that the 
Fourier-Deligne transform is the functor $\cd(E)@>>>\cd(E')$ given by 
$A\m s'_!(s^*(A)\ot\s^*\cu)[[\dim E/2]]$.

We fix a perfect symmetric bilinear pairing $\la,\ra:LG\T LG@>>>\kk$ which is invariant 
under the adjoint action of $G$. 

In this section we fix $n\in\D$. For any $\Gi$-orbit $\co$ on $L_nG$ let $\bco$ be the 
closure of $\co$ in $L_nG$. The natural $\Gi$-action on $L_nG$ has only finitely many 
orbits \cite{\GR, 3.5}. Let $\oL_nG$ be the unique open $\Gi$-orbit on $L_nG$. 

\subhead 2.2\endsubhead
Let $V$ be an algebraic variety with a given family $\ff$ of simple perverse sheaves with 
the following property: any complex in $\ff$ comes from a mixed complex on $V$ relative to
a rational structure on $V$ over some $F_q$. Let $\cd^\ff(V)$ be the subcategory of 
$\cd(V)$ whose objects are complexes $K$ such that for any $j$, any composition factor of 
${}^pH^j(K)$ is in $\ff$. Let $\ck^\ff(V)$ be the free $\ca$-module with basis $\BB^\ff(V)$
given by the isomorphism classes of simple perverse sheaves in $\ff$. Let $K$ be an object
of $\cd^\ff(V)$ with a given mixed structure relative to a rational structure of $V$ over
some $F_q$. We set
$$gr(K)=\sum_A\sum_{j,h\in\ZZ}
(-1)^j(\text{mult. of $A$ in }{}^pH^j(K)_h)(-v)^{-h}A\in\ck^\ff(V),$$
where $A$ runs over a set of representatives for the isomorphism classes in $\ff$ and the 
subscript $h$ denotes the subquotient of pure weight $h$ of a mixed perverse sheaf. (This 
agrees with the definition in \cite{\ADVIII, 36.8} after the change of variable 
$v\m(-v)\i$.)
Note that $gr(K[[m/2]])=v^m gr(K)$ for $m\in\ZZ$. 

Now let $V_1$ be another algebraic variety with a given family $\ff_1$ of simple perverse 
sheaves like $\ff$ for $V$. Then $\cd^{\ff_1}(V_1)$ is defined. Assume that we are given a 
functor $\Th:\cd(V)@>>>\cd(V_1)$ which restricts to a functor
$\cd^\ff(V)@>>>\cd^{\ff_1}(V_1)$. Assume also that $\Th$ is a composition of functors of  
the form $a_!,a^*$ induced by various maps $a$ between algebraic varieties. In particular,
$\Th$ preserves the triangulated structures and makes sense also on the mixed level. Define
an $\ca$-linear map $gr(\Th):\ck^\ff(V)@>>>\ck^{\ff_1}(V_1)$ by the following requirement:
if $A\in\ff$ is regarded as a pure complex of weight $0$ (relative to a rational structure
of $V$ over some $F_q$) then $gr(\Th)(A)=gr(\Th(A))$ where $\Th(A)$ is regarded as a mixed
complex on $V_1$ (with mixed structure defined by that of $A$). Note that $gr(\Th)(A)$ does
not depend on the choice of mixed structure. If $\Th':\cd(V_1)@>>>\cd(V_2)$ is another 
functor like $\Th$ then so is $\Th'\Th$ and we have 

(a) $gr(\Th'\Th)=gr(\Th')gr(\Th)$.

\subhead 2.3\endsubhead
For an algebraic variety $V$ with a fixed rational structure over some $F_q$ and a mixed 
complex $K$ on $V$, we define a function $\c_K:V(F_q)@>>>\bbq$ by
$$\c_K(x)=\sum_j(-1)^j(\text{trace of Frobenius on }\ch^j_x(K)).$$

\subhead 2.4\endsubhead
Let $\tfI_{L_nG}$ be the collection consisting of all irreducible $\Gi$-equivariant local
systems on various $\Gi$-orbits in $L_nG$. Let $\fI_{L_nG}$ be the set of all isomorphism
classes of irreducible $\Gi$-equivariant local systems on various $\Gi$-orbits in $L_nG$. 
For $\cl\in\tfI_{L_nG},\k\in\fI_{L_nG}$ we write $\cl\in\k$ instead of "$\k$ is the
isomorphism class of $\cl$".

For $\cl\in\k$ (as above) we say that $\cl$ or $\k$ is {\it cuspidal} (resp. {\it 
semicuspidal}) if $(G,\io)$ is rigid, $S_\k=\oL_nG$ and there exists
$\cf\in\ct_G^{cu}$ (resp. $\cf\in\ct_G^{pr}$) such that $\oL_nG\sub S_\cf$, 
$\cl\cong\cf|_{S_\k}$. 

On the other hand, if $\cf\in\ct_G^{cu}$ and $S_\cf\cap L_nG\ne\em$ then $(G,\io)$ is 
rigid and $\cf|_{\oL_nG}$ is irreducible, cuspidal in $\fI_{L_nG}$. (See 
\cite{\GR, 4.4}.)

We write $\ck(L_nG)$, $\BB(L_nG)$ instead of $\ck^\ff(L_nG)$, $\BB^\ff(L_nG)$ (see 2.2)
where $\ff$ is the family of simple $\Gi$-equivariant perverse sheaves on $L_nG$. The 
notation $\ck(L_nG)$, $\BB(L_nG)$ agrees with that in \cite{\GR, 3.9}.

For $\k\in\fI_{L_nG}$ we set $\uk^\bul=IC(\bS_\k,\cl)[[\dim S_\k/2]]$ (extended by $0$ on
$L_nG-\bS_\k$) where $\cl\in\k$. We have $\BB(L_nG)=\{\uk^\bul;\k\in\fI_{L_nG}\}$. We set
$$B(L_nG)=\{\uk;\k\in\fI_{L_nG}\}.$$

We define a $\ZZ$-linear involution $\b:\ck(L_nG)@>>>\ck(L_nG)$ by 
$\b(v^m\uk^\bul)=v^{-m}\uk^\bul$ for $m\in\ZZ$, $\k\in\fI_{L_nG}$. 

We choose a rational structure for $G$ over some $F_q$ with Frobenius map $F:G@>>>G$ such 
that $\io(t^q)=F(\io(t))$ for any $t\in\kk^*$, such that any $\Gi$-orbit in $L_nG$ or
$L_{-n}G$ is defined over $F_q$ and such that any irreducible $\Gi$-equivariant local 
system over such an orbit admits a mixed structure. Then $\Gi$ is defined over $F_q$. We
assume as we may that $\Gi$ is $F_q$-split and any connected component of $\cp^\io$ is 
defined over $F_q$. 

Let $\k\in\fI_{L_nG}$ and let $\cl\in\k$. Let $i:S_\k@>>>L_nG$ be the inclusion. We choose
a mixed structure for $\cl$ which is pure of weight $0$. Then 
$i_!\cl[[\dim S_\k/2]]$ is naturally a mixed complex on $L_nG$ and 
$$\uk:=gr(i_!\cl[[\dim S_\k/2]])\in\ck(L_nG)$$
is defined as in 2.2. It is independent of the choice of rational/mixed structures. Using 
the definitions and the purity statement in \cite{\GR, 18.2} we see that
$$\uk^\bul=\sum_{\k'\in\fI_{L_nG}}f_{\k,\k'}\uk'\tag a$$
where we have (in $\ZZ[v]$):
$$f_{\k,\k'}=\sum_{i'}(\text{mult. of $\k'$ in the local system }
\ch^iIC(\bS_\k,\cl)|_{S_{\k'}})v^{\dim S_\k-\dim S_{\k'}-i'}$$
if $S_{\k'}\sub\bS_\k$,
$$f_{\k,\k'}=0\text{ if }S_{\k'}\not\sub\bS_\k.$$
In particular, 

$f_{\k,\k}=1$,

$f_{\k,\k'}=0$ if $S_{\k'}=S_\k,\k'\ne\k$,

$f_{\k,\k'}\in v\ZZ[v]$ if $\k'\ne\k$.
\nl
We see that the $B(L_nG)$ is an $\ca$-basis of $\ck(L_nG)$.

\subhead 2.5 Induction\endsubhead
Let $P\in\cp^\io$. Then $\ck(L_n\bP)$ is defined as in 2.4 (in terms of $\bP,\io$ instead 
of $G,\io$). Now $P^\io$ and its subgroup $U_P^\io$ act freely on $\Gi\T L_nP$ by 
$y:(g,x)\m(gy\i,\Ad(y)x)$; we form the quotients $E'=\Gi\T_{U_P^\io}L_nP$, 
$E''=\Gi\T_{P^\io}L_nP$. Let $\p:L_nP@>>>L_n\bP$ be the canonical map. We have a diagram
$$L_n\bP@<a<<E'@>b>>E''@>c>>L_nG$$
where $a(g,x)=\p(x)$, $b(g,x)=(g,x)$, $c(g,x)=\Ad(g)x$.

Note that $a$ is smooth with connected fibres of dimension $s=\dim L_0P+\dim L_nU_P$, $b$ 
is a principal $\bP^\io$-bundle and $c$ is proper.

Let $A$ be a simple $\bP^\io$-equivariant perverse sheaf on $L_n\bP$. There is a well 
defined simple perverse sheaf $\tA$ on $E''$ such that
$$a^*A[[s/2]]=b^*\tA[[\dim\bP^\io/2]].$$ 
Moreover, if we regard $A$ as a pure complex of weight zero (relative to a rational 
structure over some $F_q$) then $\tA$ is naturally pure of weight zero and $c_!\tA$ is 
naturally a mixed complex whose perverse cohomology sheaves are $\Gi$-equivariant. Hence
$gr(c_!\tA)\in\ck(L_nG)$ is well defined; it is independent of the choice of mixed 
structure for $A$. Now $A\m gr(c_!\tA)$ defines an $\ca$-linear map 
$$\ind_P^G:\ck(L_n\bP)@>>>\ck(L_nG).$$

\subhead 2.6\endsubhead
Now assume that $A=IC(L_n\bP,\cl)[[\dim L_n\bP/2]]$ where $\cl\in\tfI_{L_n\bP}$ is 
cuspidal. Let 
$$\dL_nG=\{(gP^\io,z)\in\Gi/P^\io\T L_nG;\Ad(g\i)z\in\p\i(\oL_n\bP)\}.$$
We have a diagram
$$\oL_n\bP@<\ta<<\{(g,z)\in\Gi\T L_nG;\Ad(g\i)z\in\p\i(\oL_n\bP)\}@>\tb>>\dL_nG
@>\tc>>L_nG$$
with 
$$\ta(g,z)=\p(\Ad(g\i)z),\tb(g,z)=(gP^\io,z),\tc(gP^\io,z)=z.$$
Let $\dcl$ be the local system on $\dL_nG$ defined by $\ta^*\cl=\tb^*\dcl$. Using 
\cite{\GR, 4.4(b)}, we see as in \cite{\GR, 6.6} that 
$c_!\tA=\tc_!\dcl[[\dim L_0U_P/2+\dim L_nP/2]]$. If $\cl$ is regarded as a pure local 
system of weight zero (relative to a rational structure over some $F_q$) then $\dcl,A,\tA$
are naturally mixed of weight zero and in $\ck(L_nG)$ we have 
$$gr(c_!\tA)=gr(\tc_!\dcl[[\fra{\dim L_0U_P}{2}+\fra{\dim L_nP}{2}]])=
v^{\dim L_0U_P+\dim L_nP}gr(c'_!\dcl).$$

\subhead 2.7\endsubhead
We now fix $P,P'\in\cp^\io$. Let $\cl\in\tfI_{L_n\bP}$ (resp. $\cl'\in\tfI_{L_n\bP'}$) be
cuspidal. Let 
$$\align&A=IC(L_n\bP,\cl)[[\dim L_n\bP/2]]\in\cd(L_n\bP),\\&
A'=IC(L_n\bP',\cl')[[\dim L_n\bP'/2]]\in\cd(L_n\bP').\endalign$$
Let $\dL_nG,\dcl,\tc,c,\p$ be as in 2.5, 2.6, and let $\dL'_nG,\dcl',\tc',c',\p'$ be the 
analogous entities defined in terms of $P',\cl'$.

Let $R=\{h\in\Gi;hPh\i\text{ and $P'$ have a common Levi}\}$. For $h\in R$ we set 
$Q=hPh\i$; we have isomorphisms
$$\bP@>d>>\bQ@<e<<(Q\cap P')/U_{Q\cap P'}@>f>>\bP'$$
($d$ is induced by $\Ad(h)$, $e$ and $f$ are induced by the inclusions $Q\cap P'\sub Q$,
$Q\cap P'\sub P'$). Then $fe\i d:\bP@>>>\bP'$ is an isomorphism compatible with the 
homomorphisms $\io:\kk^*@>>>\bP,\io:\kk^*@>>>\bP'$. It induces a Lie algebra isomorphism 
$L\bP@>\si>>L\bP'$ compatible with the gradings hence an isomorphism
$\oL_n\bP@>\si>>\oL_n\bP'$. This carries $\cl$ to
a local system ${}^h\cl$ on $\oL_n\bP'$. We set
$$\t(h)=\dim\fra{L_nU_{P'}+\Ad(h)(L_nU_P)}{L_nU_{P'}\cap\Ad(h)(L_nU_P)}
-\dim\fra{L_0U_{P'}+\Ad(h)(L_0U_P)}{L_0U_{P'}\cap\Ad(h)(L_0U_P)}.$$
       
Let $R'$ be the set of all $h\in R$ such that ${}^h\cl\cong\cl'{}^*$. Note that $R$ and 
$R'$ are unions of $(P'{}^\io,P^\io)$-double cosets in $\Gi$ and that $\t(h)$ depends only
on the $(P'{}^\io,P^\io)$-double coset $\Om$ that contains $h$; we shall write also 
$\t_\Om$ instead of $\t(h)$.

\subhead 2.8\endsubhead
In the setup of 2.7, we choose a rational structure for $G$ over some $F_q$ with Frobenius
map $F:G@>>>G$ as in 2.4. We assume as we may that $F(P)=P,F(P')=P'$ and that 
$F^*\cl\cong\cl$, $F^*\cl'\cong\cl'$. 

For various varieties $X$ connected with $G$ which inherit an $F_q$-rational structure from
$G$ we shall write $X^F$ instead of $X(F_q)$.

We may assume that $\cl,\cl'$ have mixed structures such that all values of 
$\c_\cl:(\oL_n\bP)^F@>>>\bbq$, $\c_{\cl'}:(\oL_n\bP')^F@>>>\bbq$ are roots of $1$. Then 
$\cl,\cl',A,A'$ are pure of weight $0$. Note that 
$\dcl,\dcl',\tc_!\dcl,\tc'_!\tcl',\tA,\tA',c_!\tA,c'_!\tA'$ hence also 
$\tc_!\dcl\ot\tc'_!\tcl',c_!A\ot c'_!\tA'$ are naturally mixed complexes. We have 
$$\c_{c_!\tA}=(\sqq)^{-\dim L_0U_P-\dim L_nP}\c_{\tc_!\dcl},\qua
\c_{c'_!\tA'}=(\sqq)^{-\dim L_0U_{P'}-\dim L_nP'}\c_{\tc'_!\dcl'}.$$

\proclaim{Lemma 2.9} We have          
$$\sum_{x\in(L_nG)^F}\c_{c_!\tA\ot c'_!\tA'}(x)=
(\fra{\vt_{\Gi}}{\vt_{Z_{\bP}^0}}\sum_\Om\e_\Om(-v)^{\t_\Om})|_{v=-1/\sqq}.\tag a$$
Here $\Om$ runs over the ${P'{}^\io,P^\io}$-double cosets in $\Gi$ such that $\Om\sub R'$
and $\e_\Om$ are roots of $1$. Note that $\fra{\vt_{\Gi}}{\vt_{Z_{\bP}^0}}\in\ZZ[v^2]$.
\endproclaim

Let $N$ be the left hand side of (a). We have
$$N=(\sqq)^\vp\sum_{x\in(L_nG)^F}\c_{\tc_!\dcl\ot\tc'_!\tcl'}(x)=
(\sqq)^\vp\sum_{x'\in X^F}\c_{\dcl\bxt\dcl'}(x')$$
where 
$$\vp=-\dim L_0U_P-\dim L_nP-\dim L_0U_{P'}-\dim L_nP',$$
$$\align X=\dL_nG\T_{L_nG}\dL'_nG&
=\{(gP^\io,g'P'{}^\io,z)\in \Gi/P^\io\T\Gi/P'{}^\io\T L_nG;
\\&\Ad(g\i)z\in\p\i(\oL_n\bP),\Ad(g'{}\i)z\in\p'{}\i(\oL_n\bP')\}.\endalign$$
We have a partition $X=\sqc_\Om X_\Om$ into locally closed subvarieties indexed by the 
various $(P'{}^\io,P^\io)$-double cosets $\Om$ in $\Gi$, where 
$$\align X_\Om=&\{(gP^\io,g'P'{}^\io,z)\in\Gi/P^\io\T\Gi/P'{}^\io\T L_nG;
g'{}\i g\in\Om,\\&\Ad(g\i)z\in\p\i(\oL_n\bP),\Ad(g'{}\i)z\in\p'{}\i(\oL_n\bP')\}.\endalign
$$
(There are only finitely many such $\Om$ since $P'{}^\io,P^\io$ are parabolic subgroups of
$\Gi$.) From our assumptions we see that each $\Om$ is defined over $F_q$. We have 
$N=\sum_\Om N_\Om$ where 
$$N_\Om=(\sqq)^\vp
\sum_{(gP^\io,g'P'{}^\io,z)\in X_\Om^F}\c_{\dcl}(gP^\io,z)\c_{\dcl'}(g'P'{}^\io,z).$$
We fix an $\Om$ as above. Now $G^{\io F}$ acts on $(\dL_nG)^F$ by 
$\tg:(gP^\io,z)\m(\tg gP^\io,\Ad(\tg)z)$ and from the definitions we see that $\c_{\dcl}$ 
is constant on the orbits of this action. A similar property holds for $\c_{\dcl'}$. It 
follows that the function 
$$X_\Om^F@>>>\bbq,\qua(gP^\io,g'P'{}^\io,z)\m\c_{\dcl}(gP^\io,z)\c_{\dcl'}(g'P'{}^\io,z)$$
is constant on the orbits of the $G^{\io F}$-action on $X_\Om^F$ given by
$$\tg:(gP^\io,g'P'{}^\io,z)\m(\tg gP^\io,\tg g'P'{}^\io,\Ad(\tg)z).$$

Since $\a:X_\Om^F@>>>(\Gi/P'{}^\io)^F$, $(gP^\io,g'P'{}^\io,z)\m g'P'{}^\io$ is compatible
with the obvious actions of $G^{\io F}$ and since the $G^{\io F}$-action on 
$(\Gi/P'{}^\io)^F$ is transitive, we see that for $y\in(\Gi/P'{}^\io)^F$, the sum 
$$\sum_{(gP^\io,g'P'{}^\io,z)\in\a\i(y)}\c_{\dcl}(gP^\io,z)\c_{\dcl'}(g'P'{}^\io,z)$$
is independent of the choice of $y$. It follows that
$N_\Om=(\sqq)^\vp|(\Gi/P'{}^\io)^F|N'_\Om$ where 
$$N'_\Om=\sum\Sb(gP^\io,z)\in(\dL_nG)^F;\\ z\in\p'{}\i(\oL_n\bP),g\in\Om\eSb
\c_{\dcl}(gP^\io,z)\c_{\cl'}(\p'(z)).$$
(We have used that $\c_{\dcl'}(P'{}^\io,z)=\c_{\cl'}(\p'(z))$ which follows from the
definitions.) We set 
$$Y_\Om=\{(gP^\io,z)\in\Gi/P^\io\T L_nP';\Ad(g\i)z\in\p\i(\oL_n\bP),g\in\Om\}.$$
Define $\s:Y_\Om@>>>L_n\bP'$ by $\s(gP^\io,z)=\p'(z)$. Then
$N'_\Om=\sum_{\x\in(\oL_n\bP')^F}N''_\Om(\x)\c_{\cl'}(\x)$ where 
$$N''_\Om(\x)=\sum_{(gP^\io,z)\in\s\i(\x)^F}\c_{\dcl}(gP^\io,z).$$

Let $K_\Om=\s_!(\dcl|_{Y_\Om})$. This is naturally a mixed complex over $L_n\bP'$ and we 
have $N''_\Om(\x)=\c_{K_\Om}(\x)$ for $\x\in(\oL_n\bP')^F$. If we assume that 
$hPh\i\cap P'$ contains no Levi of $hPh\i$ for some/any $h\in\Om$, then we have $K_\Om=0$ 
(see \cite{\GR, 8.4(b)}); hence $\c_{K_\Om}=0$ and $N''_\Om(\x)=0$ for any 
$\x\in(\oL_n\bP')^F$. It follows that $N'_\Om=0$ hence $N_\Om=0$.

On the other hand, if we assume that $hPh\i\cap P'$ contains no Levi of $P'$ for some/any 
$h\in\Om$, then we have again $N_\Om=0$. (This follows from the previous paragraph applied
to $P,P',\Om\i$ instead of $P',P,\Om$.)

Assume that $\Om$ is not as in the previous two paragraphs. Thus, setting $Q=hPh\i$ for 
some $h\in\Om^F$, the intersection $Q\cap P'$ contains a Levi of $Q$ and also a Levi of 
$P'$; it follows that $Q,P'$ have a common Levi. We have $F(Q)=Q$ and $Q\in\cp^\io$. We 
have isomorphisms 
$$\bP@>d>>\bQ@<e<<(Q\cap P')/U_{Q\cap P'}@>f>>\bP'$$
($d$ is induced by $\Ad(h)$, $e$ and $f$ are induced by the inclusions $Q\cap P'\sub Q$,
$Q\cap P'\sub P'$). Then $fe\i d:\bP@>>>\bP'$ is an isomorphism compatible with the 
homomorphisms $\io:\kk^*@>>>\bP,\io:\kk^*@>>>\bP'$. It induces an isomorphism of Lie
algebras $L\bP@>\si>>L\bP'$ compatible with the gradings hence an isomorphism
$\oL_n\bP@>\si>>\oL_n\bP'$. This isomorphism carries $\cl$ to a mixed local system 
${}^h\cl$ on $\oL_n\bP'$. By \cite{\GR, 8.4(a)} and its proof \cite{\GR, 8.8}, in which 
Tate twists must be also taken in account, we see that 
$K_\Om|_{\oL_n\bP'}={}^h\cl[[-\d_\Om]]$ where 
$$\d_\Om=\dim L_0U_{P'}-\dim(L_0U_{P'}\cap L_0U_Q)+\dim(L_nU_{P'}\cap L_nU_Q)$$
and $K_\Om|_{L_n\bP'-\oL_n\bP'}=0$. We see that 
$$N'_\Om=q^{\d_\Om}\sum_{\x\in(\oL_n\bP')^F}\c_{{}^h\cl}(\x)\c_{\cl'}(\x)
=q^{\d_\Om}\sum_{\x\in(\oL_n\bP')^F}\c_{{}^h\cl\ot\cl'}(\x).\tag b$$

If ${}^h\cl\not\cong\cl'{}^*$ then ${}^h\cl\ot\cl'$ has no direct summand isomorphic to 
$\bbq$ hence by an argument as in \cite{\ADV, 23.5} we have
$H^j_c(\oL_n\bP',{}^h\cl\ot\cl')=0$ for any $j$; it follows that the last sum is $0$ so 
that $N'_\Om=0$. (To use \cite{\ADV, 23.5} we need to know that the transitive action of 
$\bP'{}^\io/Z_{\bP'}^0$ on $\oL_n\bP'$ has isotropy groups with unipotent identity 
components; in fact in our case the isotropy groups are finite as we can see from 
\cite{\GR, 4.4, 2.5(a)}.) It follows that $N_\Om=0$.

We now assume that ${}^h\cl\cong\cl'{}^*$. Then ${}^h\cl\ot\cl'$ has a unique direct
summand isomorphic to $\bbq$ and Frobenius acts on the stalk of this direct summand at any
point in $\oL_n\bP'{}^F$ as multiplication by a root of unity $\e(h)$. By an argument in 
\cite{\ADV, 24.14} we see that 
$$\sum_{\x\in\oL_n\bP'{}^F}\c_{{}^h\cl\ot\cl'}(\x)=\e(h)|\bP'{}^{\io F}||Z_{\bP'}^{0F}|\i.
$$
Hence $N_\Om=
\e(h)(\sqq)^{\vp+2\d_\Om}|(\Gi/P'{}^\io)^F||\bP'{}^{\io F}/Z_{\bP'}^{0F}|$ that
is,
$$N_\Om=\e(h)(\sqq)^{\vp+2\d_\Om-2\dim L_0U_{P'}}|G^{\io F}||Z_{\bP'}^{0F}|\i.\tag c$$

We set 

$s_0=\dim L_0U_P=\dim L_0U_Q,s'_0=\dim L_0U_{P'}$, 

$s_n=\dim L_nU_P=\dim L_nU_Q,s'_n=\dim L_nU_{P'}$,

$t'_0=\dim(L_0U_{P'}\cap L_0U_Q)$, $t'_n=\dim(L_nU_{P'}\cap L_nU_Q)$,

$r_0=\dim(L_0U_{P'}+L_0U_Q)$, $r_n=\dim(L_nU_{P'}+L_nU_Q)$.
\nl
We have $r_0+t'_0=s_0+s'_0$, $r_n+t'_n=s_n+s'_n$.

Since $Q^\io,P'{}^\io$ are parabolic subgroups of $\Gi$ with a common Levi, we have 
$s_0=s'_0$. Let $z=\dim Z_{\bP}^0=\dim Z_{\bP'}^0$. Let $g=\dim\Gi$. Since the action of 
$\bP^\io/Z_{\bP}^0$ on $L_n\bP$ has an open orbit with finite stabilizers, we have 
$\dim L_n\bP=\dim\bP^\io-z$. Moreover, $\dim\bP^\io=\dim\Gi-2\dim U_P^\io=g-2s_0$. Thus 
$\dim L_n\bP=g-z-2s_0$. Similarly, $\dim L_n\bP'=g-z-2s_0$. We have 
$\dim L_nP=\dim L_n\bP+\dim L_nU_P=g-z-2s_0+s_n$. Similarly, $\dim L_nP'=g-z-2s_0+s'_n$. We
see that the exponent of $\sqq$ in (c) is
$$\align&\vp+2\d_\Om-2\dim L_nU_{P'}\\&=-s_0-(g-z-2s_0+s_n)-s_0-(g-z-2s_0+s'_n)-2t'_0+2t'_n
\\&=2(s_0-t'_0)-(s_n+s'_n-2t'_n)-2g+2z\\&=(r_0-t'_0)-(r_n-t'_n)-2g+2z=-\t_\Om-2g+2z.\endalign$$
Thus we have
$$N_\Om=\e(h)(\sqq)^{-\t_\Om}\fra{|G^{\io F}|q^{-g}}{|Z_{\bP'}^{0F}|q^{-z}}.$$
The lemma follows.

\subhead 2.10\endsubhead
Let $\ff_2$ be the family $\{\bbq\}$ of simple perverse sheaves on a point. Write
$\ck(\pt)$ instead of $\ck^{\ff_2}(\text{point})$. We identify $\ck(\pt)=\ca$ in an obvious
way. Define an $\ca$-bilinear pairing 
$$(:):\ck(L_nG)\T\ck(L_nG)@>>>\ca$$
by the requirement that, if $K,K'$ are simple $\Gi$-equivariant perverse sheaves on $L_nG$,
we have
$$(K:K')=gr(\r_!i^*(K\bxt K'))\in\ca\tag a$$
where $K,K'$ are regarded as pure complexes of weight zero (relative to a rational 
structure over some $F_q$), $i:L_nG@>>>L_nG\T L_nG$ is the diagonal and $\r:L_nG@>>>\pt$ is
the obvious map (so that $\r_!i^*(K\bxt K')$ is a mixed complex). Note that $(K:K')$ does 
not depend on the choices hence it is well defined.

\subhead 2.11\endsubhead
Let $P,P'\in\cp^\io$. Let $V=(\Gi\T_{P^\io}L_nP)\T(\Gi\T_{P'{}^\io}L_nP')$,
$V_1=L_nG\T L_nG$, $V_2=\pt$. Let $\ff$ be the family of simple perverse sheaves on $V$ of
the form $\tA\bxt\tA'$ where $\tA$ is defined in terms of a simple $\bP^\io$-equivariant
perverse sheaf $A$ on $L_n\bP$ as in 2.5 and $\tA'$ is defined in a similar way in terms of
$A'$, a simple $\bP'{}^\io$-equivariant perverse sheaf on $L_n\bP'$. Let $\ff_1$ be the 
family of simple $\Gi\T\Gi$-equivariant perverse sheaves on $V_1$. Let $\ff_2$ be as in 
2.10. Let $\Th=(c,c')_!:\cd(V)@>>>\cd(V_1)$ where $c:\Gi\T_{P^\io}L_nP@>>>L_nG$ is as in 
2.5 and $c':\Gi\T_{P'{}^\io}L_nP'@>>>L_nG$ is the analogous map. Let 
$\Th'=\r_!i^*:\cd(V_1)@>>>\cd(V_2)$ with $\r,i$ as in 2.10. Then 2.2(a) is applicable. Thus
we have $gr(\r_!i^*(c,c')_!)=gr(\r_!i^*)gr((c,c')_!)$. If $A,A',\tA,\tA'$ are as above,
then
$$(\ind_P^G(A):\ind_{P'}^G(A'))=gr(\r_!i^*)gr((c,c')_!)(\tA\bxt\tA')$$
hence 
$$(\ind_P^G(A):\ind_{P'}^G(A'))=gr(\r_!i^*(c,c')_!)(\tA\bxt\tA').\tag a$$
We apply this to $A,A',\tA,\tA'$ as in 2.8. We choose an $F_q$-rational structure on $G$ 
and mixed structures on $A,A'$ as in 2.8. From (a) we see that
$$(\ind_P^G(A):\ind_{P'}^G(A'))=\sum_{j,h\in\ZZ}d_{j,h}(-1)^j(-v)^{-h}$$
where $d_{j,h}=\dim H^j_c(L_nG,c_!\tA\ot c'_!\tA')_h$ and the subscript $h$ denotes the 
subquotient of pure weight $h$ of a mixed vector space. Let 
$\{\l_{j,h;k};k\in[1,d_{j,h}]\}$ be the eigenvalues of the Frobenius map on 
$H^j_c(L_nG,c_!\tA\ot c'_!\tA')_h$. By the Grothendieck trace formula for the $s$-th power
of the Frobenius map ($s\in\ZZ_{>0}$), we have
$$\sum_{j,h\in\ZZ}\sum_{k\in[1,d_{j,h}]}(-1)^j\l_{j,h,k}^s=
\sum_{x\in L_nG(F_{q^s})}\c_{c_!\tA\ot c'_!\tA'}(x)$$  
(in the right hand side $\c$ is taken relative to $F_{q^s}$). Using Lemma 2.9 and its proof
(with $q$ replaced by $q^s$) we deduce
$$\sum_{j,h\in\ZZ}\sum_{k\in[1,d_{j,h}]}(-1)^j\l_{j,h,k}^s=
\sum_{\Om\in R'}(\sum_{i\in[1,u]}q^{-sa_i}-\sum_{i'\in[1,u']}q^{-sb_{i'}})
\sum_\Om\e_\Om^s(\sqq)^{-s\t_\Om}.$$
Here we write 
$$\vt_{\Gi}/\vt_{Z_{\bP}^0}=\sum\Sb i\in[1,u]\eSb v^{2a_i}
-\sum\Sb i'\in[1,u']\eSb v^{2b_{i'}}$$
with $a_1,a_2,\do,a_u,b_1,b_2,\do,b_{u'}$ in $\NN$. We can find some integer $m\ge1$ such 
that $\e_\Om^m=1$ for any $\Om\in R'$. Then for any $s\in m\ZZ_{>0}$ we have
$$\align&\sum\Sb j,h\in\ZZ\\j\text{ even}\endSb\sum_{k\in[1,d_{j,h}]}\l_{j,h,k}^s
+\sum_{\Om\in R'}\sum_{i'\in[1,u']}(\sqq)^{-2sb_{i'}-s\t_\Om}\\&
=\sum\Sb j,h\in\ZZ\\j\text{ odd}\endSb\sum_{k\in[1,d_{j,h}]}\l_{j,h,k}^s
+\sum_{\Om\in R'}\sum_{i\in[1,u]}(\sqq)^{-2sa_i-s\t_\Om}.\endalign$$
It follows that the multisets
$$\cup\Sb j,h\in\ZZ\\j\text{ even}\endSb\cup_{k\in[1,d_{j,h}]}\{\l_{j,h,k}\}
\cup\cup_{\Om\in R'}\cup_{i'\in[1,u']}\{(\sqq)^{-2b_{i'}-\t_\Om}\},$$
$$\cup\Sb j,h\in\ZZ\\j\text{ odd}\endSb\cup_{k\in[1,d_{j,h}]}\{\l_{j,h,k}\}
\cup\cup_{\Om\in R'}\cup_{i\in[1,u]}\{(\sqq)^{-2a_i-\t_\Om}\}$$
coincide. Hence for any $h\in\ZZ$ these two multisets contain the same number of elements 
$\x$ of weight $h$ (that is such that any complex absolute value of $\x$ is $(\sqq)^h$).
Since $\l_{j,h,k}$ has weight $h$, we see that for any $h$ we have
$$\align&\sum_{j\in2\ZZ}d_{j,h}+|\{(\Om,i')\in R'\T[1,u'];-2b_{i'}-\t_\Om=h\}|\\&=
\sum_{j\in2\ZZ+1}d_{j,h}+|\{(\Om,i)\in R'\T[1,u];-2a_i-\t_\Om=h\}|.\endalign$$
It follows that
$$\align&\sum_{j,h\in\ZZ;j\text{ even}}d_{j,h}(-v)^{-h}+\sum_{\Om\in R'}\sum_{i'\in[1,u']}
v^{2b_{i'}}(-v)^{\t_\Om}\\&=\sum_{j,h\in\ZZ;j\text{ odd}}d_{j,h}(-v)^{-h}
+\sum_{\Om\in R'}\sum_{i\in[1,u]}v^{2a_i}(-v)^{\t_\Om}.\endalign$$
Equivalently,
$$\sum_{j,h\in\ZZ}d_{j,h}(-1)^j(-v)^{-h}=
\fra{\vt_{\Gi}}{\vt_{Z_{\bP}^0}}\sum_{\Om\in R'}(-v)^{\t_\Om}$$
that is,
$$(\ind_P^G(A):\ind_{P'}^G(A'))=
\fra{\vt_{\Gi}}{\vt_{Z_{\bP}^0}}\sum_{\Om\in R'}(-v)^{\t_\Om}.\tag b$$

\subhead 2.12\endsubhead
The pefect pairing $\la,\ra:LG\T LG@>>>\kk$ (see 2.1) restricts to a perfect pairing 
$L_nG\T L_{-n}G@>>>\kk$ denoted again by $\la,\ra$. Note that 
$\la\Ad(g)x,\Ad(g)x'\ra=\la x,x'\ra$ for $x\in L_nG,x'\in L_{-n}G,g\in\Gi$. The
Fourier-Deligne transform $\cd(L_nG)@>>>\cd(L_{-n}G)$ (defined as in 2.1 in terms of 
$\la,\ra$) takes a simple $\Gi$-equivariant perverse sheaf $A$ on $L_nG$ to a simple 
$\Gi$-equivariant perverse sheaf $\Ph^G_n(A)$ (or $\Ph_n(A)$) on $L_{-n}G$. Moreover 
$A\m\Ph_n(A)$ defines a bijection $\BB(L_nG)@>\si>>\BB(L_{-n}G)$ and this extends uniquely
to an isomorphism $\Ph_n:\ck(L_nG)@>\si>>\ck(L_{-n}G)$ of $\ca$-modules. From
\cite{\GR, 3.14(a)} we see that this is the inverse of the isomorphism 
$\Ph_{-n}:\ck(L_{-n}G)@>\si>>\ck(L_nG)$ defined like $\Ph_n$ in terms of $-n$ instead of 
$n$. Let 
$$\k\m\dk,\qua\fI_{L_nG}@>>>\fI_{L_{-n}G}$$
be the bijection such that for any $\k\in\fI_{L_nG}$ we have 
$$\Ph_n(\uk^\bul)=\un{\dk}^\bul.\tag a$$
The inverse of this bijection is denoted again by $\k\m\dk$.

For any simple $\Gi$-equivariant perverse sheaf $A$ on $L_nG$, the restriction of the
$\Gi$-action on $LG_n$ to $Z_G$ (a subgroup of $\Gi$) is the trivial action of $Z_G$. 
Then $Z_G$ acts naturally by automorphisms on $A$ and this action is via scalar 
multiplication by a character $c_A:Z_G@>>>\kk^*$ (trivial on the identity component of 
$Z_G$). From the definitions we see that 

(b) $c_{\Ph_n(A)}=c_A$.
\nl
Now assume that $(G,\io)$ is rigid and that $\cl\in\tfI_{L_nG}$ is cuspidal (so that
$S_\cl=\oL_nG$). Let $A=IC(L_nG,\cl)[[\dim L_nG/2]]$. According to \cite{\GR, 10.6(e)} we 
have $\Ph_n(A)=IC(L_{-n}G,\cl')[[\dim L_{-n}G/2]]$ where $\cl'\in\tfI_{L_{-n}G}$ is 
cuspidal (so that $S_\cl=\oL_{-n}G$). We can find $\cf\in\ct_G^{cu},\cf'\in\ct_G^{cu}$ such
that $\cl=\cf|_{\oL_nG}$, $\cl'=\cf'|_{\oL_{-n}G}$. We show: 

(c) $\cf=\cf'$.
\nl
From (b) we see that the natural action of $Z_G$ on any stalk of $\cl$ and the
natural action of $Z_G$ on any stalk of $\cl'$ is through the same character of $Z_G$.
Now $Z_G$ also acts naturally on each stalk of $\cf$ and on each stalk of $\cf'$ though
some character of $Z_G$ (in the adjoint action of $G$ on $LG$, $Z_G$ acts trivially). Since
$\cl=\cf|_{\oL_nG}$, $\cl'=\cf'|_{\oL_{-n}G}$, the character of $Z_G$ attached to $\cf$ is
the same as the character of $Z_G$ attached to $\cf'$. But from the classification of
cuspidal local systems it is known that an object in $\ct_G^{cu}$ is completely determined 
by the associated character of $Z_G$. This proves (c).

\subhead 2.13\endsubhead
Let $P\in\cp^\io$. Let $\wInd_{LP}^{LG}:\ck(L_n\bP)@>>>\ck(L_nG)$ be the $\ca$-linear map
defined in \cite{\GR, 6.2}. We show:
$$\wInd_{LP}^{LG}(\x)=\ind_P^G(\x)\text{ for any }\x\in\ck(L_n\bP).\tag a$$
Let $A$ be a simple $\bP^\io$-equivariant perverse sheaf on $L_n\bP$. Let $\tA$ and 
$c:E''@>>>L_nG$ be as in 2.5. We regard $A$ as a pure complex of weight zero (relative to a
rational structure over some $F_q$). Then $\tA$ is naturally pure of weight zero and, by
Deligne \cite{\DE}, $c_!\tA$ is pure of weight zero. Using \cite{\BBD} we deduce that for 
any $j$, ${}^pH^j(c_!\tA)$ is pure of weight $j$. Hence the definition of $\ind_P^G(A)$ 
becomes
$$\ind_P^G(A)=\sum_{A_1}\sum_{j\in\ZZ}(\text{mult. of $A_1$ in }{}^pH^j(c_!\tA))v^{-j}A_1
\tag b$$
where $A_1$ runs over the set of simple $\Gi$-equivariant perverse sheaves on $L_nG$ (up to
isomorphism). On the other hand, since $c_!\tA$ has weight zero, we have (by \cite{\BBD})
$c_!\tA\cong\op_{j\in\ZZ}{}^pH^j(c_!\tA)[-j]$ in $\cd(L_nG)$ hence
$\wInd_{LP}^{LG}(A)$ is equal to the right hand side of (b). We see that 
$\ind_P^G(A)=\wInd_{LP}^{LG}(A)$. This proves (a).

From (a) we see that a number of results proved in \cite{\GR} for $\wInd_{LP}^{LG}$ imply
corresponding results for $\ind_P^G$ (see (c)-(f) below).

(c) {\it The elements $\ind_P^G(A)\in\ck(L_nG)$ with $P\in\cp^\io$ and $A$ as in 2.6 span
the $\QQ(v)$-vector space $\QQ(v)\ot_\ca\ck(L_nG)$.}
\nl
(We use \cite{\GR, 13.3, 17.3}.)

(d) {\it Let $P,P'\in\cp^\io$. Assume that $P\sub P'$. Let $Q$ be the image of $P$ under 
$P'@>>>\bP'$. Note that $Q$ is a parabolic subgroup of $\bP'$ containing $\io(\kk^*)$
and $\bQ=\bP$.) Then $\ind_P^G:\ck(L_n\bP)@>>>\ck(L_nG)$ is equal to the composition}

$\ck(L_n\bP)@>\ind_Q^{\bP'}>>\ck(L_n\bP')@>\ind_{P'}^G>>\ck(L_nG)$.
\nl
(We use \cite{\GR, 6.4}.)

(e) {\it Let $P\in\cp^\io$. For any $\x\in\ck(L_n\bP)$ we have 
$\ind_P^G(\Ph^{\bP}_n(\x))=\Ph^G_n(\ind_P^G(\x))\in\ck(L_{-n}G)$ (here $\ind_P^G$ in the
left hand side is defined in terms of $-n$ instead of $n$).}
\nl
(We use \cite{\GR, 10.5}.)

(f) {\it Let $P\in\cp^\io$, $\x\in\ck(L_n\bP)$. We have 
$\ind_P^G(\b(\x))=\b(\ind_P^G(\x))$.}
\nl
The proof is along the lines of the proof of the analogous equality 
\cite{\ADVIII, 36.9(c)} which is based on the relative hard Lefschetz theorem \cite{\BBD}.

\subhead 2.14\endsubhead
Let $\ZZ((v))$ be the ring of power series $\sum_{j\in\ZZ}a_jv^j$ ($a_j\in\ZZ$) such that 
$a_j=0$ for $j\ll 0$. We have naturally $\ca\sub\ZZ((v))$ and $\ZZ((v))$ becomes an 
$\ca$-algebra. 

Let $\{,\}:\ck(L_nG)\T\ck(L_nG)@>>>\ZZ((v))$ be the $\ca$-bilinear pairing defined in
\cite{\GR, 3.11}. We show:
$$(\x:\x')=\vt_{\Gi}\{\x,\x'\}|_{v\m-v}\text{ for any }\x,\x'\in\ck(L_nG).\tag a$$
If $\x=\ind_P^G(A),\x'=\ind_{P'}^G(A')$ with $P,P',A,A'$ as in 2.7, then (a) follows by
comparing 2.11(b) with the analogous formula for 
$\{\wInd_{LP}^{LG}(A),\wInd_{LP'}^{LG}(A')\}$ in \cite{\GR, 15.3}. (We use also 2.13(a).) 
This implies (a) in general, in view of 2.13(c).

For $\k\in\fI_{L_nG},\cl\in\k$, we define $\k^*\in\fI_{L_nG}$ by $\cl^*\in\k^*$. Then
$\un{\k^*}^\bul$ is the Verdier dual of $\uk^\bul$. We show that for $\k,\k'$ in 
$\fI_{L_nG}$ we have
$$(\uk^\bul:\uk'{}^\bul)\in\d_{\k',\k^*}+v\ZZ[v].\tag b$$
$$(\uk:\uk')\in\d_{\k',\k^*}+v\ZZ[v].\tag c$$
Using (a) and \cite{\GR, 3.11(d)}, we see that to prove (b) it is enough to verify the 
following statement. 

If $f\in\ca$ and $\vt_{\Gi}\i f\in\d+v\ZZ[[v]]$ with $\d\in\ZZ$ then $f\in\d+v\ZZ[v]$.
\nl
This is clear since $\vt_{\Gi}\i\in1+v\ZZ[[v]]$.

Now (c) follows from (b) using the fact that the transition matrix from $(\uk)$ to
$(\uk^\bul)$ is uni-triangular with off-diagonal entries in $v\ZZ[v]$ (see 2.4(a)).

\subhead 2.15\endsubhead
Let $x\in L_nG$. Let $\co$ be the $\Gi$-orbit of $x$ in $L_nG$.  Let $P=P(x)$ be the 
parabolic subgroup associated to $x$ in \cite{\GR, 5.2}. Recall that $\io(\kk^*)\sub P$. We
show:

(a) {\it The adjoint action of $U_P^\io$ on $x+L_nU_P$ is transitive.}
\nl
Let $S$ be the orbit of $x$ under this action. Since $S$ is an orbit of an action of a 
unipotent group on the affine space $x+L_nU_P$, it is closed in $x+L_nU_P$. Hence it is
enough to show that $\dim S=\dim L_nU_P$ or that 
$\dim U_P^\io-\dim(U_P^\io)_x=\dim L_nU_P$, where $(U_P^\io)_x$ is the stabilizer of $x$ in
$U_P^\io$. The last equality is proved in the course of the proof of \cite{\GR, 5.9}. This
proves (a).

Let $\p:L_nP@>>>L_n\bP$ be the canonical map. We have $x\in\p\i(\oL_n\bP)$, see
\cite{\GR, 5.3(b)}. We show:

(b) {\it The adjoint action of $P^\io$ on $\p\i(\oL_n\bP)$ is transitive.}
\nl
Let $y\in\p\i(\oL_n\bP)$. To show that $y$ is in the $P^\io$-orbit of $x$ we may replace 
$y$ by a $P^\io$-conjugate. Hence we may assume that $y\in x+L_nU_P$. In that case we may 
use (a). This proves (b).

Let $L_n\bP@<a<<E'@>b>>E''@>c>>L_nG$ be as in 2.5 (defined in terms of the present $P$).
Let $E''_1=\Gi\T_{P^\io}\p\i(\oL_n\bP)$, an open subset of $E''$. We show:

(c) {\it $c:E''@>>>L_nG$ restricts to an isomorphism $E''_1@>>>\co$.}
\nl 
Using (b) we see that $\Gi$ acts transitively on $E''_1$ hence $c(E''_1)$ is a single 
$\Gi$-orbit. It contains $x$ hence it equals $\co$. We see that $E''_1\sub c\i(\co)$. By
the proof of \cite{\GR, 6.8(b)}, $c$ restricts to an isomorphism $c\i(\co)@>\si>>\co$. In
particular, $c\i(\co)$ is a single $\Gi$-orbit. Since $E''_1$ is a $\Gi$-orbit contained in
$c\i(\co)$ we must have $E''_1=c\i(\co)$ and (c) follows.

Now let $\cl'$ be an irreducible $\bP^\io$-equivariant local system on $\oL_n\bP$. Define
$\k'\in\fI_{L_n\bP}$ by $\cl'\in\k'$. Then $\uk'\in\ck(L_n\bP)$ is defined as in 2.4. Let 
$\tcl'$ be the local system on $E''_1$ whose inverse image under the obvious map
$\Gi\T_{U_P^\io}\p\i(\oL_n\bP)@>>>E''_1$ coincides with the inverse image of $\cl'$ under

$\Gi\T_{U_P^\io}\p\i(\oL_n\bP)@>>>\oL_n\bP_n$, $(g,x)\m\p(x)$.
\nl
Let $\cl$ be the irreducible 
$\Gi$-equivariant local system on $\co$ corresponding to $\tcl'$ under the isomorphism 
$E''_1@>>>\co$ induced by $c$. (See (c).) Define $\k\in\fI_{L_nG}$ by $\cl\in\k$. Then 
$\uk\in\ck(L_nG)$ is defined as in 2.4. We show:

(d) $\uk=\ind_P^G(\uk')$.
\nl
We choose an $F_q$-rational structure on $G$ as in 2.4 so that $x$ is $F_q$-rational hence
$P$ is defined over $F_q$, and a mixed structure for $\cl'$ which is pure of weight $0$. 
Let $i':\oL_n\bP_n@>>>L_n\bP$, $i:\co@>>>L_nG$, $i_1:E''_1@>>>E''$ be the inclusions. Let 
$(A_r)_{r\in[1,m]}$ be a set of representatives for the simple $\bP^\io$-equivariant
perverse sheaves on $L_n\bP$. Now $A=i'_!\cl'[[\dim\oL_n\bP/2]])$ is naturally a mixed 
complex on $L_n\bP$ and we have 
$$\uk'=\sum_{j,h}\sum_{r\in[1,m]}(-1)^j(\text{mult. of $A_r$ in }{}^pH_j(A)_h)(-v)^{-h}A_r
\in\ck(L_n\bP).$$
We attach to each $A_r$ a simple perverse sheaf $\tA_r$ on $E''$ by

$a^*A_r[[s/2]]=b^*\tA_r[[\dim\bP^\io/2]]$
\nl
where $s$ is as in 2.5. Let $\tA=i_{1!}\tcl'[[\dim E''_1/2]]\in\cd(E'')$. From the 
definitions we have $a^*A[[s/2]]=b^*\tA[[\dim\bP^\io/2]]$. Since $a,b$ are smooth morphisms
with connected fibres of dimension $s,\dim\bP^\io$ we deduce that
$$a^*({}^pH^j(A)_h)[[s/2]]=b^*({}^pH^j(\tA)_h)[[\dim\bP^\io/2]]$$
and 
$$\text{mult. of $A_r$ in }{}^pH_j(A)_h=\text{mult. of $\tA_r$ in }{}^pH_j(\tA)_h$$
for any $r,j,h$. Hence
$$\uk'
=\sum_{j,h}\sum_{r\in[1,m]}(-1)^j(\text{mult. of $\tA_r$ in }{}^pH_j(\tA)_h)(-v)^{-h}A_r.$$
By the definition of $\ind_P^G$ we have (in $\ck(L_nG)$):
$$\ind_P^G(\uk')=\sum_{j,h}\sum_{r\in[1,m]}(-1)^j
(\text{mult. of $\tA_r$ in }{}^pH_j(\tA)_h)(-v)^{-h}gr(c_!(\tA_r))=gr(c_!\x)$$
where 
$$\x=\sum_{j,h}\sum_{r\in[1,m]}
(-1)^j(\text{mult. of $\tA_r$ in }{}^pH_j(\tA)_h)(-v)^{-h}\tA_r=gr(\tA)\in\ck^\ff(E'')$$
and $\ff$ is the family of simple perverse sheaves on $E''$ of the form $\tA_r(r\in[1,m])$.
Thus,
$$\ind_P^G(\uk')=gr(c_!(gr(i_{1!}\tcl'[[\dim E''_1/2]]))).$$
Let $\ff_1$ be the family of simple $\Gi$-equivariant perverse sheaves on $E''_1$. Let 
$\ff_0$ be the family of simple $\Gi$-equivariant perverse sheaves on $L_nG$. Applying 
2.2(a) to $\Th=i_{1!}:\cd^{\ff_1}(E''_1)@>>>\cd^\ff(E'')$, 
$\Th'=c_!:\cd^\ff(E'')@>>>\cd^{\ff_0}(L_nG)$, we obtain $gr(c_!i_{1!})=gr(c_!)gr(i_{1!})$.
We see that
$$\ind_P^G(\uk')=gr(c_!i_{1!}\tcl'[[\dim\co/2]]).$$
(Recall that $\dim E''_1=\dim\co$.) From the definitions we have $c_!i_{1!}\tcl'=i_!\cl$. 
Hence 
$$\ind_P^G(\uk')=gr(i_!\cl[[\dim\co/2]])=\uk.$$
This proves (d).

\subhead 2.16\endsubhead
We preserve the setup of 2.15. Let $\co_2$ be a $\bP^\io$-orbit in $L_n\bP-\oL_n\bP$. Let 
$E''_2=\Gi\T_{P^\io}\p\i(\co_2)$, a subset of $E''$. We show:

(a) {\it The image of $E''_2$ under $c:E''@>>>L_nG$ is contained in $\bco-\co$.}
\nl
Let $y\in c(E'')$. We show that $y\in\bco$. We have $y=\Ad(g)\et$ for some 
$g\in\Gi$ and $\et\in L_nP$. Replacing $y$ by $\Ad(g\i)y$ we may assume that 
$y\in L_nP$. By \cite{\GR, 5.9}, $L_nP$ is contained in the closure of the $P^\io$-orbit of
$x$ in $L_nP$ which is clearly contained in $\bco$. Thus $y\in\bco$. We see that
$c(E'')\sub\bco$. In particular, $c(E''_2)\sub\bco$. By the proof of 2.15(c) we 
have $E''_1=c\i(\co)$. Since $E''_1\cap E''_2=\em$, we have $c\i(\co)\cap E''_2=\em$ hence
$c(E''_2)\cap\co=\em$. Thus, $c(E''_2)\sub\bco-\co$ and (a) is proved.

Now let $\k''\in\fI_{L_n\bP}$ be such that $S_{\k''}=\co_2$. Then $\uk''\in\ck(L_n\bP)$ is
defined as in 2.4. From the definitions we see that $\ind_P^G(\uk'')\in\sum_\k\ca\uk$ where
$\k$ runs over the elements of $\fI_{L_nG}$ such that $S_\k$ is contained in the closure of
$c(E''_2)$. Using this and (a) we see that

(b) $\ind_P^G(\uk'')\in\sum_{\k;S_\k\sub\bco-\co}\ca\uk$.

\subhead 2.17\endsubhead
Let $\cl',\cl''\in\tfI_{L_nG}$. Define $\k',\k''\in\fI_{L_nG}$ by $\cl'\in\k'$,
$\cl''\in\k''$. Assume that $S_{\cl'}\cap S_{\cl''}=\em$. We show:
$$(\uk':\uk'')=0.\tag a$$
We choose an $F_q$-rational structure on $G$ as in 2.4 and mixed structures on $\cl',\cl''$
which makes them pure of weight $0$. Let $i':S_{\cl'}@>>>L_nG$, $i'':S_{\cl''}@>>>L_nG$.
Let $V=S_{\cl'}\T S_{\cl''},V_1=L_nG\T L_nG,V_2=\pt$. Let $\ff_1,\ff_2$ be as in 2.11. Let
$\ff$ be the family of simple perverse sheaves on $V$ consisting of $\cl'\bxt\cl''$. Let
$\Th=(i',i'')_!:\cd(V)@>>>\cd(V_1)$. Let $\Th'=\r_!i^*:\cd(V_1)@>>>\cd(V_2)$ where $\r,i$
are as in 2.10. Then 2.2(a) is applicable. Thus we have 
$gr(\r_!i^*(i',i'')_!)=gr(\r_!i^*)gr((i',i''')_!)$. Hence
$$gr(\r_!i^*)gr((i',i''')_!)(\cl'\bxt\cl'')=gr(\r_!i^*(i',i'')_!)(\cl'\bxt\cl'').$$
The right hand side is zero since $i^*(i',i'')_!(\cl'\bxt\cl'')=0$ (by our assumption
$S_{\cl'}\cap S_{\cl''}=\em$). Thus, $gr(\r_!i^*)gr((i',i''')_!)(\cl'\bxt\cl'')=0$. Using 
this and 2.10(a) we see that (a) holds.

\head 3. Computation of multiplicities\endhead
\subhead 3.1\endsubhead
For any $n\in\D$ let $\Gi\bsl L_nG$ be the (finite) set of $\Gi$-orbits on $L_nG$. We 
define a map 

(a) $\fP_n@>>>\Gi\bsl L_nG$
\nl
as follows. Let $P\in\fP_n$. Let $M,L^r_tG$ be as in 1.9. We have $L_0^0=LG'$ for a well 
defined connected reductive subgroup $G'$ of $G$. Now $G'$ acts on $L^n_nG$ (by restriction
of the adjoint action of $G$ on $LG$) and there is a unique open $G'$-orbit $\co_0$ for 
this action. Since $G'\sub\Gi$, there is a unique $\Gi$-orbit $\co$ on $L_nG$ that contains
$\co_0$. Now (a) associates $\co$ to $P$. Clearly (a) factors through a map 

(b) $\ufP_n@>>>\Gi\bsl L_nG$.
\nl
This is a bijection. Its inverse associates to the $\Gi$-orbit of $x\in L_nG$ the
$\Gi$-orbit of the parabolic subgroup associated to $x$ in \cite{\GR, 5.2}. 

Another parametrization of $\Gi\bsl L_nG$ was given by Vinberg \cite{\VI} (see also 
\lb Kawanaka \cite{\KA}).

For any $\et\in\ufP_n$ let $\co_\et$ be the $\Gi$-orbit in $L_nG$ corresponding to 
$\et$ under (b). Note that

(c) $\dim\co_\et=d_\et$
\nl
where $d_\et$ is as in 1.10. This follows easily from \cite{\GR, 5.4(a), 5.9}.

\subhead 3.2\endsubhead
In the remainder of this section we assume that $\kk$ is as in 2.1. From the bijection 
3.1(b) we see that $B(L_nG)=\sqc_{\et\in\ufP_n}B^\et_n$ where
$$B^\et_n:=\{\uk;\k\in\fI_{L_nG},S_\k=\co_\et\}.$$

\subhead 3.3\endsubhead
We set ${}^{\QQ(v)}\ck(L_nG)=\QQ(v)\ot_\ca\ck(L_nG)$.
For $n\in\D$ we define a $\QQ(v)$-linear map
$$t_n=t_n^G:K_G@>>>{}^{\QQ(v)}\ck(L_nG)$$
by sending the basis element $\II_\cs$ to $\ind_P^G(A)$ where $(P,\ce)\in\cs$ and 

$A=IC(L_n\bP,\ce|_{\oL_n\bP})[[\dim L_n\bP/2]]$. 

The pairing ${}^{\QQ(v)}\ck(L_nG)\T{}^{\QQ(v)}\ck(L_nG)@>>>\QQ(v)$ obtained from the
pairing $(:)$ in 2.10 by linear extension will be denoted again by $(:)$. From 2.11(b) we 
see that the equality 

(a) $(t_n(\x):t_n(\x'))=(\x:\x')$
\nl
holds when $\x,\x'$ run through a basis of $K_G$; hence it holds for any $\x,\x'$ in $K_G$.
Since $t_n$ is surjective (see 2.13(c)) and the pairing $(:)$ on ${}^{\QQ(v)}\ck(L_nG)$ is
non-degenerate (see 2.14(b)) we see that $\ker t_n=\car_G$ so that $t_n$ induces an 
isomorphism
$$\tit_n:K_G/\car_G@>>>{}^{\QQ(v)}\ck(L_nG).\tag b$$ 

\subhead 3.4\endsubhead
Let $n\in\D$. We extend $\b:\ck(L_nG)@>>>\ck(L_nG)$ to a $\QQ$-linear endomorphism of
${}^{\QQ(v)}\ck(L_nG)$ (denoted again by $\b$) by $\r\ot\x\m\bar\r\ot\b(\x)$. We show that

(a) $\b(t_n(\x))=t_n(\b(\x))$
\nl
for any $\x\in K_G$. We may assume that $\x=\II_\cs$ for some $\cs\in\ucj_G$. Then we have
$\b(\x)=\x$. It is enough to show that $\b(\ind_P^G(A))=\ind_P^G(A)$ where $P,A$ are as in
3.3. By 2.13(f) we have $\b(\ind_P^G(A))=\ind_P^G(\b(A))$ and it remains to note that
$\b(A)=A$.

We extend $\Ph_n:\ck(L_nG)@>>>\ck(L_{-n}G)$ to a $\QQ(v)$-linear isomorphism \lb
${}^{\QQ(v)}\ck(L_nG)@>>>{}^{\QQ(v)}\ck(L_{-n}G)$ (denoted again by $\Ph_n$).
We show that 

(b) $\Ph_n(t_n(\x))=t_{-n}(\x)$
\nl
for any $\x\in K_G$. We may assume that $\x=\II_\cs$ for some $\cs\in\ucj_G$. It is then 
enough to show that if $(P,\ce)\in\cs$ and 

$A=IC(L_n\bP,\ce|_{\oL_n\bP})[[\dim L_n\bP/2]]$,
$A'=IC(L_{-n}\bP,\ce|_{\oL_{-n}\bP})[[\dim L_{-n}\bP/2]]$,
\nl
then 
$\Ph^G_n(\ind_P^G(A))=\ind_P^G(A')$. Using 2.13(e) we see that it is enough to show that 
$\ind_P^G(\Ph^{\bP}_n(A))=\ind_P^G(A')$. Hence it is enough to show that
$\Ph^{\bP}_n(A)=A'$. This follows from 2.12(c) (applied to $\bP$ instead of $G$).

\subhead 3.5\endsubhead
Let $n\in\D$. Let $Q\in\cp^\io$. For any $\x\in K_{\bQ}$ we have

(a) $t_n^G(f_Q^G(\x))=\ind_Q^G(t_n^{\bQ}(\x))$.
\nl
where $f_Q^G$ is as in 1.4 and $\ind_Q^G$ is extended by $\QQ(v)$-linearity.
We may assume that $\x=\II_{\cs'}$ for some 
$\cs'\in\ucj_{\bQ}$; then the result follows from 2.13(d).

We show that 1.8(c) holds. Using 3.3(b) we see that it is enough to show that
$t_n(f_Q^G(\car_{\bQ}))=0$. Using (a) it is enough to show that
$\ind_Q^G(t_n^{\bQ}(\car_{\bQ}))=0$. This follows from $t_n^{\bQ}(\car_{\bQ})=0$ (see 3.3).

\subhead 3.6\endsubhead
In Section 1 we tried to associate to any $n\in\D$ and $\et\in\ufP_n$ a subset 
$\cz^\et_n$ of $K_G$. We will go again through the definitions (with the help of results 
in Section 2) and we will add the requirement that

(a) {\it for any $\et\in\ufP_n$, $t_n$ restricts to a bijection} $\cz^\et_n@>\si>>B^\et_n$.
\nl
We may assume that $G$ is not a torus and that the subsets $\cz^\et_n$ are already defined
when $G$ is replaced by any $\bP$ with $P\in\fP'_n$. (If $G$ is a torus then $\et$
must be $\{G\}$ and we define $\cz^\et_n$ to be the subset consisting of the unique basis 
element of $K_G$.)

Assume first that $n\in\D$ and $\et\in\ufP'_n$ (see 1.9). We define $\cz^\et_n$ as in 1.11.
We show that (a) holds for our $\et$. Let $\x\in\cz^\et_n$. With notation in 1.11(i) we 
have $\x=f_P^G(\x')$ for some $\x'\in\cz^{\{\bP\}}_n$ where $P\in\et$. Using 3.5(a),
$t_n^G(f_P^G(\x'))$ is equal to $\ind_P^G(t_n^{\bP}(\x'))$ which by the induction 
hypothesis belongs to $\ind_P^G(B^{\{\bP\}}_n)$. 
Thus, $t_n(f_P^G(\x'))\in\ind_P^G(B^{\{\bP\}}_n)$.
By 2.15(d) and its proof, $\ind_P^G$ maps $B^{\{\bP\}}_n$ 
into $B^\et_n$ and in fact defines a
bijection $a:B^{\{\bP\}}_n@>\si>>B^\et_n$ 
(we use the isomorphism in \cite{\GR, 5.8}). Thus we
have $t_n(f_P^G(\x'))\in B^\et_n$. We consider the diagram
$$\CD
\cz^{\{\bP\}}_n@>a_1>>\cz^\et_n\\
@Va_2VV                   @Va_3VV    \\
B^{\{\bP\}}_n@>a>>B^\et_n \endCD$$
where $a_1$ defined by $\ind_P^G$, $a_2$ is defined by $t_n^{\bP}$ and $a_3$ is defined by 
$t_n^G$. This diagram is commutative by 3.5(a). By the induction hypothesis, $a_2$ is a 
bijection. We have just seen that $a$ is a bijection. It follows that $a_3a_1$ is a 
bijection. Hence $a_1$ is injective. By the definition of $\cz^\et_n$, $a_1$ is surjective.
Thus, $a_1$ is a bijection. In particular, 1.12(a) holds.
Since $a_3a_1$ is a bijection we see that $a_3$ is a bijection.
Thus (a) holds in our case.

\subhead 3.7\endsubhead
Let $n\in\D$. Define $\cz'_n$ as in 1.11. Let $L'_nG=\{x\in L_nG;P(x)\ne G\}$. Let 
$$I_n=\{\k\in\fI_{L_nG};S_\k\sub L'_nG\},\qua \uI_n=\{\uk;\k\in I_n\}.$$
Now $t_n$ defines a bijection $\cz'_n@>\si>>\uI_n$. 

We show that 1.12(b) holds. Let $\et,\et'$ be two distinct elements of $\ufP'_n$. Then 
$\cz^\et_n,\cz^{\et'}_n$ are 
disjoint since their images $B^\et_n$, $B^{\et'}_n$ under $t_n$ are disjoint. (A local
system in $B^\et_n$ has a support different from that of a local system in $B^{\et'}_n$
since 3.1(b) is a bijection.) 

We show that 1.12(c) holds. Using 3.5(a) and 3.6(a) we see that this follows from 2.16(b).
(We use also 3.1(c).)

We show that 1.12(d) holds. Using 3.4(a) and 3.6(a) we see that it is enough to prove the 
following statement (for $G$ instead of $\bP$). If $(G,\io)$ is rigid and
$\k\in\fI_{L_nG}-I_n$ then $\b(\k)-\k\in\sum_{\k'\in I_n}\ca\k'$. This is immediate from
the definitions.

We show that 1.12(e) holds. Using 3.3(a) we see that it is enough to show that the matrix 
with entries $(t_n(\x):t_n(\x'))$ indexed by $\cz'_n\T\cz'_n$ is non-singular. It is also
enough to show that the matrix with entries $(\uk:\uk')$ indexed by $I_n\T I_n$ is 
non-singular. This follows from 2.14(c) since $I_n$ is stable under $\k\m\k^*$.

\subhead 3.8\endsubhead
Let $[\uI_n]$ be the $\ca$-submodule of $\ck(L_nG)$ with basis $\uI_n$. Now $L_nG-L'_nG$ is
empty (resp. is $\oL_nG$) if $(G,\io)$ is not rigid (resp. rigid). Hence $L'_nG$ is a 
closed subset of $L_nG$. This, together with 2.4(a) shows that 
$\{\uk^\bul;\k\in I_n\}$ is an $\ca$-basis of $[\uI_n]$. 

\subhead 3.9\endsubhead
For $\x\in\cz'_n$ we define $W_n^\x$ as in 1.13. We have $t_n(\x)=\uk$ where $\k\in I_n$. 
We show:

(a) $t_n(W_n^\x)=\uk^\bul$.
\nl
Let $y=t_n(W_n^\x)$. Applying $t_n$ to the equality $\b(W_n^\x)=W_n^\x\mod\car_G$ in 1.13 
and using 3.4(a) we see that $\b(y)=y$. Applying $t_n$ to the equality 
$W_n^\x=\sum_{\x_1\in\cz'_n}c_{\x,\x_1}\x_1$ in 1.13 we obtain
$y=\sum_{\x_1\in\cz'_n}c_{\x,\x_1}t_n(\x_1)$. We see that $y=\sum_{\k'}\tf_{\k,\k'}\uk'$ 
where $\k'$ runs over the elements in $\fI_{L_nG}$ and

$\tf_{\k,\k'}\ne0$ implies $\dim S_{\k'}<\dim S_\k$ or $\x=\x'$;

$\tf_{\k,\k'}\ne0$, $\k\ne\k'$ implies $\tf_{\k,\k'}\in v\ZZ[v]$;

$\tf_{\k,\k'}=1$ if $\k=\k'$.
\nl
These conditions together with $\b(y)=y$ determine $y$ uniquely.  Since $\uk^\bul$ 
satisfies the same conditions as $y$ (see 2.4) we see that $y=\uk^\bul$. This proves (a).

\subhead 3.10\endsubhead
Until the end of 3.11 we assume that $(G,\io)$ is rigid. For $n\in\D$ we define  
$J_{-n}$ as in 1.15. Let $C=\QQ(v)\ot_\ca[\uI_n]$. Let $\x_0\in\cz'_{-n}$. Then 
$t_{-n}(\x_0)=\uk$ where $\k\in I_{-n}$. We have $\x_0\in J_{-n}$ if and only if 
$t_n(W_{-n}^{\x_0})\n C$ that is, if and only if $\Ph_n(t_n(W_{-n}^{\x_0}))\n\Ph_n(C)$ that
is (using 3.4(b)), if and only if $t_{-n}(W_{-n}^{\x_0}))\n\Ph_n(C)$ that is (using 3.9(a) 
with $n$ replaced by $-n$), if and only if $\uk^\bul\n\Ph_n(C)$. Now 
$\{\uk_1^\bul;\k_1\in I_n\}$ is a $\QQ(v)$-basis of $C$. By 2.12(a), 
$\{\un{\dk_1}^\bul;\k_1\in I_n\}$ is a $\QQ(v)$-basis of $\Ph_n(C)$. Hence the 
condition that $\uk^\bul\in\Ph_n(C)$ is equivalent to the condition that $\k=\dk_1$ for
some $\k\in I_n$. We see that $\x_0\in J_{-n}$ if and only if $\dk\n I_n$.

Assssume now that $\x_0\in J_{-n}$. Since $\dk\n I_n$, we have
$S_{\dk}=\oL_nG$. Define $h_n$ as in 1.15. Let 
$z=t_nh_n(\x_0)$. From the definitions we have $t_n(W_{-n}^{\x_0})-z\in C$ and $(z:C)=0$. 
As we have seen earlier we have $\Ph_n(t_n(W_{-n}^{\x_0}))=\uk^\bul$. Hence
$t_n(W_{-n}^{\x_0})=\Ph_{-n}(\uk^\bul)=\un{\dk}^\bul$. Thus we have
$\un{\dk}^\bul-z\in C$. Since $S_{\dk}=\oL_nG$ we have
$\un{\dk}^\bul=\un{\dk}\mod C$. It follows that $\un{\dk}-z\in C$. Using 
$S_{\dk}=\oL_nG$ and 2.17(a) we see that $(\un{\dk}:C)=0$. Since $(z:C)=0$ we see 
that $(\un{\dk}-z:C)=0$. Since $(:)$ is non-degenerate on $C$ (see 3.7) and
$\un{\dk}-z\in C$ we see that $\un{\dk}-z=0$. Thus 

(a) {\it $t_nh_n$ is the map $\x_0\m\un{\dk}$ where $\k\in I_{-n}$ is given by}
$t_{-n}(\x_0)=\uk$.

We show that 1.16 holds. It is enough to show that $h_n$ is injective. It is also enough to
show that $t_nh_n$ is injective. Let $\x'_0\in J_{-n}$. Define $\k'\in I_{-n}$ by
$t_{-n}(\x'_0)=\uk'$. Assume that $\un{\dk}=\un{\dk'}$. Then $\dk=\dk'$ and
$\k=\k'$. Since $t_{-n}:\cz'_{-n}@>>>I_{-n}$ is bijective it follows that $\x_0=\x'_0$. 
Thus 1.16 is proved.

Let $\cc_n=h_n(J_{-n})$ (see 1.15). The previous proof shows that the map 
$\cc_n@>>>t_n(\cc_n)$ (restriction of $t_n$) is a bijection. We see that

$t_n(\cc_n)=\{\uk';\k'\in\fI_{L_nG},S_{\k'}=\oL_nG,\dk'\in I_{-n}\}$.

\subhead 3.11\endsubhead
Define $\cc'$ as in 1.14. If $\cf\in\ct_G^{pr},S_\cf=C_G^\io$ and $r_\cf,\cs_\cf$ are as in
1.14, then for $n\in\D$ we have $t_n(r_\cf\i\II_{\cs_\cf})=\uk=\uk^\bul$ where
$\k\in\ck(L_nG)$ is $\cf|_{\oL_nG}$. (The last two equalities follow from 
\cite{\GR, 11.13}.) Replacing $n$ by $-n$ we have similarly
$t_{-n}(r_\cf\i\II_{\cs_\cf})=\uk'=\uk'{}^\bul$ where $\k'\in\ck(L_{-n}G)$ is 
$\cf|_{\oL_{-n}G}$. Using 3.4(b) we have 
$$\Ph_n(\uk^\bul)=\Ph_n(t_n(r_\cf\i\II_{\cs_\cf}))=
t_{-n}(r_\cf\i\II_{\cs_\cf}))=\uk'{}^\bul.$$
Thus $\k'=\dk$ so that $\dk\n I_{-n}$.

We show that the map $\cc'@>>>t_n(\cc')$ (restriction of $t_n$) is a bijection. It is 
enough to note that the map $\cf\m\cf|_{\oL_nG}$ is a bijection from $\ct_G^{pr}$ to the
set of semicuspidal objects in $\fI_{L_nG}$. This follows from the fact that, if 
$x\in\oL_nG$, the centralizer of $x$ in $G$ and the centralizer of $x$ in $\Gi$ have the 
same group of components. 

We now show that 1.17 holds. First we show that $\cc_n\cap\cc'=\em$. It is enough to show
that if $\uk\in t_n(\cc_n)$ and $\un{\tik}\in t_n(\cc')$ then $\k\ne\tik$. From our 
assumption we have $\dk\in I_{-n}$, $\dot{\tik}\n I_{-n}$ (see above). Thus, 
$\k\ne\tik$, as required.

Next we show that $\cz'_n\cap(\cc_n\cup\cc')=\em$. It is enough to show that if 
$\uk\in t_n(\cz'_n)$ and $\un{\tik}\in t_n(\cc_n\cup\cc')$ then $\k\ne\tik$. From our 
assumption we have $S_{\uk}\ne\oL_nG$, $S_{\un{\tik}}=\oL_nG$, Thus, $\k\ne\tik$, as 
required. This proves 1.17. We see also that the map 
$\cz'_n\cup\cc_n\cup\cc'@>>>t_n(\cz'_n\cup\cc_n\cup\cc')$ (restriction of $t_n$) is a 
bijection. 

\subhead 3.12\endsubhead
If $(G,\io)$ is not rigid then the definition of the subsets $\cz_n^\et$ 
($\et\in\ufP_n$) is complete. If $(G,\io)$ is rigid then $\ufP_n-\ufP'_n=\{G\}$. For 
$n\in\D$ we set $\cz_n^{\{G\}}=\cc_n\cup\cc'$. By 1.17, this union is disjoint. The 
definition of the subsets $\cz_n^\et$ ($\et\in\ufP_n$) is complete. Define $\cz_n$ as in 
1.18. Note that the map $\cz_n@>>>t_n(\cz_n)$ (restriction of $t_n$) is a bijection. 

We show that $t_n(\cz_n)=B(L_nG)$. Let $\k\in\fI_{L_nG}$. If $\k\in I_n$ then 
$\uk\in t_n(\cz'_n)$. If $\k\in\fI_{L_nG}-I_n$ and $\dk\in I_{-n}$ then 
$\uk\in t_n(\cc_n)$. If $\k\in\fI_{L_nG}-I_n$ and $\dk\in\fI_{L_{-n}G}-I_{-n}$ then by 
\cite{\GR, 12.3}, we have $\uk\in t_n(\cc')$. We see that 

(a) {\it $t_n$ restricts to a bijection $\cz_n@>\si>>B(L_nG)$.}

\subhead 3.13\endsubhead
For $n\in\D$ and $\x\in\cz_n$ we define an element $W_n^\x$ as in 1.19. We show:

(a) {\it If $t_n(\x)=\uk$ with $\k\in\fI_{L_nG}$ then $t_n(W_n^\x)=\uk^\bul$.}
\nl
When $\x\in\cz'_n$ this follows from 3.9(a). 
Next assume that $\x\in\cc_n$. Define $\x_0\in J_{-n}$ by $h_n(\x_0)=\x$. Define
$\k_0\in\fI_{L_{-n}G}$ by $t_{-n}(\x_0)=\uk_0$. Using the definition, 3.4(b) and 3.9(a) 
(for $-n$ instead of $n$) we have 

$t_n(W_n^\x)=t_n(W_{-n}^{\x_0})=\Ph_{-n}(t_{-n}(W_{-n}^{\x_0}))=\Ph_{-n}\uk_0^\bul=
\un{\dk_0}^\bul$.
\nl
By 3.10(a) we have $\uk=t_n(\x)=\un{\dk_0}$. Thus (a) holds in our case.

Finally, assume that $\x\in\cc'$. In this case we have 
$t_n(W_n^\x)=t_n(\x)=\uk=\uk^\bul$, see 3.11. This proves (a).

Let $\k\in\fI_{L_nG}$. Let $\x\in\cz_n$ be such that $t_n(\x)=\uk$. Let $c_{\x,\x'}$ be as
in 1.19 ($\x'\in\cz_n$). Applying $t_n$ to both sides of 1.19(a) and using (a) we obtain

$\uk^\bul=\sum_{\x'\in\cz_n}c_{\x,\x'}t_n(\x')$.
\nl
Comparing this with 2.4(a) we see that for any $\k,\k'$ in $\fI_{L_nG}$ we have
$$f_{\k,\k'}=c_{\x,\x'}\tag b$$
where $\x,\x'\in\cz_n$ are defined by $t_n\x=\k,t_n\x'=\k'$. Note that (b) provides a 
method to compute explicitly the matrix of multiplicities $(f_{\k,\k'})$.

\subhead 3.14\endsubhead
We prove 1.20. 
Let $\cs\in\ucj^D$. Let $(P,\ce)\in\cs$, $\cl=\ce|_{\oL_n\bP}$. Let $\tc,\dcl,A$ be as in 
2.6. We regard $\dcl$ as a pure local system of weight $0$. By \cite{\GR, 21.1(b)} and its
proof, $\ch^i\tc_!\dcl$ is pure of weight $i$ and is $0$ unless $i\in 2\NN$. It follows
that in $\ck(L_nG)$ we have $\ind_P^G(A)=\sum_{\k\in\fI_{L_nG}}\te_{\cs,\k}\uk$ where 
$\te_{\cs,\k}\in\ca$ is equal to a power of $v$ times
$\sum_i(\text{mult. of $\cl'$ in $\ch^i\tc_!\dcl$})v^{-i}$. (Here $\cl'\in\k$.) From the 
definitions we have for any $\x\in\cz_n$

(a) $e_{\cs,\x}=\te_{\cs,\k}$
\nl
where $\k=t_n(\x)$. Hence 1.20 holds. 

\head 4. Further results\endhead
\subhead 4.1\endsubhead
In this section we assume that $\kk$ is as in 2.1.

To any $n\in\D$ and any $\k\in\fI_{L_nG}$ we shall associate a $\Gi$-orbit 
$\cs_\k\in\ucj_G$, an element $r_\k\in\ca-\{0\}$ and an element $L_\k\in\ck(L_nG)$ such 
that $L_\k=t_n(r_\k\i\II_{\cs_\k})$. We may assume that these objects are already defined
when $G$ is replaced by $\bP$ with $P\in\fP'_n$.

(i) Assume first that $\k$ is semicuspidal. There is a unique $\cf\in\cj_G^{pr}$ such that
$S_\cf\cap L_n\G=\oL_nG$ and $\cf|_{\oL_nG}\in\k$. Let $\cs_k=\cs_\cf$ (see 1.18), 
$r_\k=r_\cf$ (see 1.3), $L_\k=\uk=\uk^\bul$. These elements satisfy the required 
condition, see 3.11.

(ii) Next we assume that $\k\in I_n$. We have $\uk\in B^\et_n$ for a unique $\et\in\ufP'_n$.
Let $P\in\et$. Then $P\ne G$. Let $a:B^{\{\bP\}}_n@>\si>>B^\et_n$ be the bijection in 3.6.
Let $\uk_1=a\i(\uk)$. Now $\cs_{\k_1}$, $r_{\k_1}$, $L_{\k_1}\in\ck(L_n\bP)$ are already 
defined from the induction hypothesis. Let $r_\k=r_{\k_1}$, $L_k=\ind_P^G(L_{\k_1})$, 
$\cs_\k=a_P^G(\cs_{\k_1})$. These elements satisfy the required condition.

(iii) Next we assume that $\k\n I_n$ and $\k$ is not semicuspidal. By 3.12 we have
$\dk\in I_{-n}$. Now $\cs_{\dk}$, $r_{\dk}$, $L_{\dk}\in\ck(L_{-n}G)$ are 
defined as in (ii). Let $\cs_\k=\cs_{\dk}$, $r_\k=r_{\dk}$, 
$L_\k=\Ph_{-n}(L_{\dk})$. These elements satisfy the required condition.

This completes the definition of $\cs_\k,r_\k,L_\k$.

\subhead 4.2\endsubhead
For $n\in\D$, we shall define a partial order $\le$ on $\fI_{L_nG}$. We may assume that 
$\le$ is already defined when $G$ is replaced by $\bP$ with $P\in\fP'_n$.

(i) Assume that at least one of $\k,\k'$ is semicuspidal. Then $\k\le\k'$ if and only if 
$\k=\k'$.

(ii) Assume that $S_\k\ne S_{\k'}$ and neither $\k$ nor $\k'$ is semicuspidal. Then 
$\k\le\k'$ if and only if $S_\k\sub\bS_{\k'}-S_{\k'}$.

(iii) Assume that $S_\k=S_{\k'}$ and $\k\in I_n$ (hence also $\k'\in I_n$). We have 
$\uk,\uk'\in B^\et_n$ for a unique $\et\in\ufP'_n$. Let $P\in\et$; then $P\ne G$. Let 
$a:B^{\{\bP\}}_n@>\si>>B^\et_n$ be the bijection in 3.6. Let 
$\uk_1=a\i(\uk),\uk'_1=a\i(\uk')$. We say that $\k\le\k'$ if and only if $\k_1\le\k'_1$
(which is known by the inductive assumption applied to $\bP$).

(iv) Assume that $S_\k=S_{\k'}$, $\k\n I_n$ (hence $\k'\n I_n$) and neither $\k$ nor $\k'$
is semicuspidal. By 3.12 we have $\dk\in I_{-n}$, $\dk'\in I_{-n}$. We say that 
$\k\le\k'$ if and only if $\dk\le\dk'$ which is known from (ii) (if 
$S_{\dk}\ne S_{\dk'}$) or (iii) (if $S_{\dk}=S_{\dk'}$).

This completes the definition of $\le$. We write $\k'<\k$ instead of $\k'\le\k,\k'\ne\k$. 

\subhead 4.3. Example\endsubhead
In this subsection we assume that $G$ is the group of automorphisms of a $4$-dimensional 
$\kk$-vector space $V$ preserving a fixed non-degenerate symplectic form. We fix a direct
sum decomposition $V=V_{-1}\op V_1$ where $V_1,V_{-1}$ are Lagrangian subspaces. For any
$t\in\kk^*$ define $\io(t)\in G$ by $\io(t)x=tx$ for $x\in V_1$, $\io(t)x=t\i x$ for 
$x\in V_{-1}$. Let $\D=\{2,-2\}$. Note that $(G,\io)$ is rigid. We identify $\Gi$ with
$GL(V_1)$ by $g\m g|_{V_1}$. The grading of $LG$ defined by $\io$ has non-zero components 
in degrees $-2,0,2$ and $L_2G$ (resp. $L_{-2}G$) may be identified as a representation of 
$\Gi$ with $S^2V_1$ (resp. $S^2V_1^*$) where $S^2$ stands for the second symmetric power. 
For $n\in\D$, the set $\fI_{L_nG}$ consists of four objects 
$\k_{0,n},\k_{2,n},\k_{3,n},\tik_{3,n}$ where $\k_{i,n}$ represents the local system $\bbq$
on the $\Gi$-orbit of dimension $i (i=0,2,3)$ and $\tik_{3,n}$ represents a non-trivial 
local system of rank $1$ on the open orbit. The effect of the Fourier-Deligne transform is
as follows.

$\Ph_n(\k_{0,n}^\bul)=\k_{3,-n}^\bul$, $\Ph_n(\k_{2,n}^\bul)=\tik_{3,-n}^\bul$, 
$\Ph_n(\k_{3,n}^\bul)=\k_{0,-n}^\bul$, $\Ph_n(\tik_{3,n}^\bul)=\k_{2,-n}^\bul$.
\nl
The partial order in 4.2 is: 

$\k_{0,n}<\k_{2,n}<\k_{3,n}<\tik_{3,n}$.
\nl
We have

$L_{\k_{0,n}}=\k_{0,n}^\bul$, $L_{\k_{2,n}}=\k_{2,n}^\bul+\k_{0,n}^\bul$,
$L_{\k_{3,n}}=\k_{3,n}^\bul$, $L_{\tik_{3,n}}=\tik_{3,n}^\bul+\k_{3,n}^\bul$.

\subhead 4.4\endsubhead
We show that for any $n\in\D$ and any $\k\in\fI_{L_nG}$ we have

(a) $L_\k\in\uk+\sum_{\k'\in\fI_{L_nG};\k'<\k}\text{$\ca$ }\uk'$.
\nl
We may assume that (a) is already known when $G$ is replaced by $\bP$ with $P\in\fP'_n$.

(i) Assume that $\k$ is semicuspidal. Then $L_\k=\k$ and (a) is clear.

(ii) Assume that $\k\in I_n$. Let $P,a,\k_1$ be as in 4.1(ii). By the induction hypothesis
we have $L_{\k_1}\in\uk_1+\sum_{\k'_1;\k'_1<\k_1}\text{$\ca$ }\uk'_1$.
Applying $\ind_P^G$ and using 2.15(d) we obtain

$L_\k\in\uk+\sum_{\k'_1\in X}\text{$\ca$ }\un{a(\k'_1)}
+\sum_{\k'_1\in Y}\text{$\ca$ }\ind_P^G(\uk'_1)$
\nl
where

$X=\{\k'_1;\k'_1<\k_1,S_{\k'_1}=S_{\k_1}\}$,
$Y=\{\k'_1;\k'_1<\k_1,S_{\k'_1}\ne S_{\k_1}\}$.
\nl
For $\k'_1\in X$ we have $S_{a(\k'_1)}=S_\k$ and $a(k'_1)<\k$. For $\k'_1\in Y$ we have 
$S_{\k'_1}\sub S_{\k_1}-S_{\k_1}$. By 2.16(b), for any $\k'_1\in Y$, $\ind_P^G(\uk'_1)$ is
an $\ca$-linear combination of elements $\uk'$ where $S_{\k'}\sub S_\k-S_\k$ (hence 
$\k'<\k$). Hence $L_\k$ satisfies (a).

(iii) Assume that $\k\n I_n$ and $\k$ is not semicuspidal. By 3.12 we have
$\dk\in I_{-n}$. By (ii), we have

$L_{\dk}\in\un{\dk}+\sum_{\k';\k'<\dk}\text{$\ca$ }\uk'$.
\nl
Using $\uk_0\in\uk_0^\bul+\sum_{\k'_0<\k_0}\text{$\ca$ }\uk'_0{}^\bul$ for 
$\k_0\in I_{-n}$ we deduce 
$L_{\dk}\in\un{\dk}^\bul+\sum_{\k';\k'<\dk}\text{$\ca$ }\uk'{}^\bul$. Applying 
$\Ph^G_{-n}$ we obtain

$L_\k\in\uk^\bul+\sum_{\k';\k'<\dk}\text{$\ca$ }\un{\dk'}^\bul$.
\nl
Using $\uk_1^\bul\in\uk_1+\sum_{\k'_1<\k_1}\text{$\ca$ }\uk'_1$ for $\k_1\in I_n$ we 
see that it is
enough to show that for any $\k'$ such that $\k'<\dk$ we have $\dk'<\k$. If 
$S_\k=S_{\dk'}$ then this follows from 4.2(iv). If $S_\k\ne S_{\dk'}$ then
$S_{\dk'}\sub S_\k-S_\k$ (we have $S_\k=\oL_nG$) hence again $\dk'<\k$ (using 4.2(ii)).

This completes the proof of (a).

From (a) we deduce

(b) {\it The set $\{L_\k;\k\in\fI_{L_nG}\}$ is an $\ca$-basis of $\ck(L_nG)$.}

\subhead 4.5\endsubhead
We show that:

(a) {\it the map $\fI_{L_nG}@>>>\ucj_G,\k\m\cs_\k$ is injective.}
\nl
Let $\k,\k'\in\fI_{L_nG}$ be such that $\cs_\k=\cs_{\k'}$. Then 
$t_n(\cs_\k)=t_n(\cs_{\k'})$ hence $r_\k L_\k=r_{\k'}L_{\k'}$. Using now 4.4(b)
we deduce that $\k=\k'$ as desired. 

\subhead 4.6\endsubhead
Let $\x\in\cz_n$ and let $\k=t_n(\x)$. By \cite{\CUIII, 3.36}, the numbers 
$\te_{\cs,\k}|_{v=1}$ (for various $\cs$, see 3.14) are the dimensions of the various 
"weight spaces" of a standard module over an affine Hecke algebra. Using 3.14(a) we see 
that the dimensions of these weight spaces are given by the numbers $e_{\cs,\x}|_{v=1}$ 
(for various $\cs$) hence are computable from the algorithm in Section 1.

\subhead 4.7\endsubhead
In this subsection we shall summarize some of the results of this paper in terms of the
vector space $\bK_G=K_G/\car_G$. (Any text marked as $\sp\do\sp$ applies only in the case
where $(G,\io)$ is rigid.) Note that the pairing $(:)$ on $K_G$ induces a pairing 
$\bK_G\T\bK_G@>>>\QQ(v)$ denoted again by $(:)$. Also $\b:K_G@>>>K_G$ induces an involution 
$\bK_G@>>>\bK_G$ denoted again by $\b$. Let $\p:K_G@>>>\bK_G$ be the obvious map.

Since for $n\in\D$, $B(L_nG)$ is a $\QQ(v)$-basis of ${}^{\QQ(v)}\ck(L_nG)$, we see (using
3.12(a) and 3.3(b)) that $\p$ restricts to a bijection of $\cz_n$ onto a basis $\bcz_n$ of 
$\bK_G$. We say that $\bcz_n$ is a {\it PBW-basis} of $\bK_G$. $\sp$ The last bijection 
restricts to bijections of $\cz'_n,\cc_n,\cc'$ onto subsets $\bcz'_n,\occ_n,\occ'$ of 
$\bcz_n$. $\sp$ In the case where $(G,\io)$ is not rigid we set $\bcz'_n=\bcz_n$.

Let $\cm_n$ be the $\ZZ[v]$-submodule of $\bK_G$ with basis $\bcz_n$. For any $u\in\bcz_n$
let $\bW^u_n=\p(W^\x_n)$ where $\x\in\cz_n$ is given by $\p(\x)=u$. From 3.13(a) and 
2.4(a) we see that $\{\bW^n_u;u\in\bcz_n\}$ is a $\ZZ[v]$-basis of $\cm_n$ and that for 
any $u\in\bcz_n$ we have $\bW^u_n-u\in v\cm_n$. Define a bijection $u\m\du$ of $\bcz_n$ 
onto $\bcz_{-n}$ as follows. Let $\x\in\cz_n$ be such that $\p(\x)=u$; let 
$\k\in\fI_{L_nG}$ be such that $t_n(\x)=\uk$ (see 3.12(a)). Let $\x'\in\cz_{-n}$ be such 
that $t_{-n}\x'=\un{\dk}$ (see 2.12). Then $\du=\p(\x')$. The
inverse of the bijection $u\m\du$ is denoted again by $u\m\du$.

For $u,\x,\k,\x'$ as above we have $\bW^u_n=\p(W^\x_n), t_n(W^\x_n)=\uk^\bul$,
$\bW^{\du}_{-n}=\p(W^{\x'}_{-n})$, $t_{-n}(W^{\x'}_{-n})=\un{\dk}^\bul$. By 2.12(a) we
have $\Ph_n(\uk^\bul)=\un{\dk}^\bul$ hence $t_{-n}(W^{\x'}_{-n})=\Ph_n(t_n(W^\x_n))$
and this equals $t_{-n}(W^\x_n))$ (see 3.4(b)). Thus, $t_{-n}(W^{\x'}_{-n}-W^\x_n)=0$. 
Since $\ker t_{-n}=\car_G$ we see that $W^{\x'}_{-n}-W^\x_n\in\car_G$. Applying $\p$ we 
deduce

(a) $\bW^{\du}_{-n}=\bW^u_n$.
\nl
Moreover, from the proof in 3.10, we see that 

(b) $\sp$ $u\in\occ_n\imp\du\in\bcz_{-n}$. $\sp$
\nl
From (a) we see that the basis $(\bW^u_n)$ of $\bK_G$ coincides with the basis 
$(\bW^u_{-n})$. We call this the {\it canonical basis} of $\bK_G$. It follows that 
$\cm_n=\cm_{-n}$. We shall write $\cm$ instead of $\cm_n=\cm_{-n}$. We have 

(c) $\bW^u_n-u\in v\cm$, $\bW^{\du}_{-n}-\du\in v\cm$. 
\nl
Combining with (a) we see that 

(d) $u-\du\in v\cm$ for any $u\in\bcz_n$.
\nl
Let $\p':\cm@>>>\cm/v\cm$ be the obvious map. From (a),(c) we see that
there exists a $\ZZ$-basis $X$ of $\cm/v\cm$ such that $\p'$ restricts to bijections 
$(\bW^u_n)=(\bW^u_{-n})@>\si>>X$, $\bcz_n@>\si>>X$, $\bcz_{-n}@>\si>>X$.
$\sp$ Moreover, $X$ can be partitioned as $X=X_0\sqc X_n\sqc X_{-n}\sqc X'$ so that $\p'$
restricts to bijections $\occ_n@>\si>>X_n$, $\occ_{-n}@>\si>>X_{-n}$, $\occ'@>\si>>X'$,
$\bcz'_n@>\si>>X_0\cup X_{-n}$, $\bcz'_{-n}@>\si>>X_0\cup X_n$. $\sp$ We set 

(e) $\tX=\bcz'_n\cup\bcz'_{-n}$ (if $(G,\io)$ is not rigid),
$\tX=\bcz'_n\cup\bcz'_{-n}\cup\occ'$ (if $(G,\io)$ is rigid).
\nl
We have

(f) $X=\p'(\tX)$.
\nl
We show that 

(g) $\tX$ generates the $\ZZ[v]$-module $\cm$.
\nl
If $(G,\io)$ is not rigid, this is clear. $\sp$ Assume now that $(G,\io)$ is rigid. Let 
$\cm'$ be the $\ZZ[v]$-submodule of $\cm$ generated by $\tX$. If $u\in\bcz'_n$ then by the
arguments in 1.13, we have $\bW^u_n\in\sum_{u'\in\bcz'_n}\ZZ[v]u'$ hence $\bW^u_n\in\cm'$.
If $u\in\occ_n$ then $\bW^u_n=\bW^{\du}_{-n}$ and this is in $\cm'$ since 
$\du\in\bcz'_{-n}$ (we use the previous sentence with $u,n$ replaced by $\du,-n$). If
$u\in\occ'$ then $\bW^u_n=u$ is again in $\cm'$. Since $(\bW^u_n)_{u\in\bcz_n}$ is a 
$\ZZ[v]$-basis of $\cm$ we see that $\cm=\cm'$. $\sp$ This proves (g).

We show how the canonical basis and the PBW-bases are determined in terms of the
subsets $\bcz'_n,\bcz'_{-n}$ (which are defined by the inductive construction in 1.11(i)) 
and (in the rigid case) by the set $\occ'$ which is defined as in 1.14.

We first define $\tX$ as in (e). Note that $\cm$ is defined in terms of $\tX$ as in (g) 
and then the basis $X$ of $\cm/v\cm$ is defined in terms of $\tX$ as in (f).

Now the canonical basis can be reconstructed in terms of $\cm$ and $X$: for any 
$x\in X$ there is a unique element $\hx\in\cm$ such that $\p'(\hx)=x$ and $\b(\hx)=\hx$. 
The elements $\{\hx;x\in X\}$ form the canonical basis. Now let $n\in\D$. We show how to 
reconstruct the PBW-basis $\bcz_n$. If $(G,\io)$ is not rigid then $\bcz_n=\bcz'_n$ is
already known. $\sp$ Assume now that $(G,\io)$ is rigid. Then the part $\bcz'_n\cup\occ'$ 
of $\bcz_n$ is already known. It remains to characterize the subset $\occ_n$ of $\bcz_n$. 
For $n\in\D$ let $X_n$ be the set of all $x\in\p'(\bcz'_{-n})$ such that 
$x\n\p'(\bcz'_n)$. For any $x\in X_n$ we can write uniquely $\hx=x'+x''$ where $x''$
is in the subspace of $\bK_G$ spanned by $\bcz'_n$ and $x'$ is orthogonal under $(:)$ to 
that subspace. Then $\occ_n$ consists of the elements $x'$ for various $x\in X_n$. $\sp$

\subhead 4.8\endsubhead
Let $n\in\D$. Assume that $(G,\io)$ is rigid. 
Let $\Xi_n$ be the set of all $\k\in\fI_{L_nG}$ such that $S_\k\ne\oL_nG$ and
$S_{\dk}=\oL_{-n}G$. Here $\k\m\dk$ is as in 2.12.
It would be interesting to find a simple description of the set of local systems
$\Xi_n$ (without using Fourier-Deligne transform). In particular, we would like
to know which $\Gi$-orbits in $L_nG$ are of the form $S_\k$ for some $\k\in\Xi_n$.
(Our results answer this question only in terms of an algorithm, not in closed 
form.) For example, if $G=GL_n(\kk)$ then $\Xi_n$ has only one object: the local 
system $\bbq$ on the $0$-dimensional orbit. In the case studied in 4.3, $\Xi_n$ 
has two objects: the local system $\bbq$ on the $0$ or $2$ dimensional orbit.

\enddocument